\documentclass{article}
\usepackage{amsfonts}
\usepackage{amsmath}
\usepackage[backref,colorlinks=true]{hyperref}
\newtheorem{theorem}{Theorem}

\newtheorem{corollary}[theorem]{Corollary}

\newtheorem{definition}[theorem]{Definition}
\newtheorem{example}[theorem]{Example}

\newtheorem{lemma}[theorem]{Lemma}

\begin{document}
\date{}
\title{ Multitype branching processes with immigration in random environment
and polling systems }
\author{ Vatutin Vladimir\thanks{Supported in part be the RFBR grant
08-01-00078 and the program "Mathematical control theory" of RAS}\\
Steklov Mathematical Institute,\\
 Gubkin street 8, 119991,\\
Moscow, Russia\\
E-mail: \tt vatutin@mi.ras.ru} \maketitle

\begin{abstract}
For multitype branching processes with immigration evolving in a
random environment and producing a final product we find the tail
distribution of the size of the final product accumulated in the
system for a life period. Using this result  we investigate the
tail distribution of the busy periods of the branching type
polling systems with random service disciplines and random
positive switch-over times.
\end{abstract}

\section{Introduction}

Branching processes with and without immigration are powerful tools in
studying various models of queueing systems (see, for instance, \cite{ViSe06}%
,\cite{Gri88},\cite{Gri90},\cite{Gri92},\cite{Mei07}, \cite{Res93}, \cite%
{VD2002} and \cite{YY07}). In this paper we use multitype
branching processes with accumulation of a final product and
immigration  which evolve in random environment (MBPFPIRE) to
study the tail distribution of busy periods of a class of polling
systems in which input parameters, service disciplines and the
distributions of switch-over times vary in a random manner. This
article complements the  results of paper \cite{Vat2009}
established for polling systems with zero switch-over times.

The paper is organized as follows. In Sections \ref{DFF}-\ref{MM}
we established various results related to MBPFPIRE. These results
 are applied in Section \ref{Spoll}
 to investigate  the tail distribution and moments of the branching type
  polling systems  with positive switch-over times.

\section{Branching processes in random environment with final
product}\label{DFF}

Let $\mathbf{s:}=(s_{1},\ldots ,s_{m})\in \left[ 0,1\right] ^{m}$ be a $m$%
-dimensional variable,
\begin{equation*}
\mathbf{s}^{\mathbf{k}}:=s_{1}^{k_{1}}\cdots s_{m}^{k_{m}},\,k_{i}\in
\mathbb{N}_{0}=\left\{ 0,1,2,\ldots \right\} ,
\end{equation*}%
and $\left( \mathbf{\xi };\varphi \right) $ be a $(m+1)$-dimensional vector,
where the vector $\mathbf{\xi :=}\left( \xi _{1},\ldots ,\xi _{m}\right) $
has integer-valued nonnegative random variables as components and $\varphi $
is a nonnegative random variable, and let%
\begin{equation*}
F(\mathbf{s;}\lambda ):=\mathbf{E}\left[ \mathbf{s}^{\mathbf{\xi }%
}e^{-\lambda \varphi }\right] =\mathbf{E}\left[ s_{1}^{\xi _{1}}s_{2}^{\xi
_{2}}\cdots s_{m}^{\xi _{m}}e^{-\lambda \varphi }\right] ,\,\mathbf{s}\in %
\left[ 0,1\right] ^{m},\,\lambda \geq 0,
\end{equation*}%
be the respective mixed probability generating function (m.p.g.f.). Denote $%
\ \mathcal{F}:=\left\{ F(\mathbf{s;}\lambda )\right\} $ the set of all such
m.p.g.f.'s and let%
\begin{equation*}
\mathcal{F}^{m}:=\left\{ \mathbf{F}(\mathbf{s;}\lambda )=\left( F^{(1)}(%
\mathbf{s;}\lambda ),\ldots ,F^{(m)}(\mathbf{s;}\lambda )\right) :F^{(i)}(%
\mathbf{s;}\lambda )\in \mathcal{F},\,i=1,...,m\right\}
\end{equation*}%
be the $m$-times direct product of $\mathcal{F}$. Let, further,
\begin{equation*}
\mathcal{F}_{0}:=\left\{ f(\mathbf{s})=F(\mathbf{s;}0)=\mathbf{E}\left[
s_{1}^{\xi _{1}}s_{2}^{\xi _{2}}\cdots s_{m}^{\xi _{m}}\right] :F(\mathbf{s;}%
\lambda )\in \mathcal{F}\right\}
\end{equation*}%
be the set of all ordinary probability generating functions (p.g.f.'s) and
\begin{equation*}
\mathcal{F}_{0}^{m}:=\left\{ \mathbf{f}(\mathbf{s})=\left( f^{(1)}(\mathbf{s}%
),\ldots ,f^{(m)}(\mathbf{s})\right) :f^{(i)}(\mathbf{s})\in \mathcal{F}%
_{0},\,i=1,...,m\right\}
\end{equation*}%
be the set of all $m$-dimensional (vector-valued) p.g.f.'s.

Let, further, $\left( \mathbf{\eta };\psi \right) $ be a $(m+1)$-dimensional
tuple where the vector $\mathbf{\eta }:=\left( \eta _{1},\ldots ,\eta
_{m}\right) $ has integer-valued nonnegative random components and $\psi $
is a nonnegative random variable and let%
\begin{equation*}
G(\mathbf{s};\lambda )=\mathbf{E}\left[\,\mathbf{s}^{\mathbf{\eta }%
}e^{-\lambda \psi }\right] =\mathbf{E}\left[\, s_{1}^{\eta
_{1}}s_{2}^{\eta
_{2}}\cdots s_{m}^{\eta _{m}}e^{-\lambda \psi }\right] ,\,\mathbf{s}\in %
\left[ 0,1\right] ^{m},\,\lambda \geq 0,
\end{equation*}%
be the respective m.p.g.f. Denote $\ \mathcal{G}:=\left\{ G(\mathbf{s;}%
\lambda )\right\} $ the set of all such m.p.g.f.'s (which, of course, is
equivalent to $\ \mathcal{F}$). Let, further,%
\begin{equation*}
\mathcal{G}_{0}:=\left\{
g(\mathbf{s})=G(\mathbf{s;}0)=\mathbf{E}\left[\,
s_{1}^{\eta _{1}}s_{2}^{\eta _{2}}\cdots s_{m}^{\eta _{m}}\right] ,\,G(%
\mathbf{s;}\lambda )\in \mathcal{G}\right\}
\end{equation*}%
be the set of all ordinary probability generating functions (p.g.f.'s) .

Assume that a probability measure $\mathbb{Q}$ is specified on the natural $%
\sigma $-algebra $\mathcal{A}$ generated by the subsets of $\mathcal{F}%
^{m}\times \mathcal{G}.$ Let $\mathbf{H}(\mathbf{s;}\lambda ),\mathbf{H}_{0}(%
\mathbf{s;}\lambda ),\mathbf{H}_{1}(\mathbf{s;}\lambda ),\ldots ,\mathbf{H}%
_{k}(\mathbf{s;}\lambda ),\ldots $ where%
\begin{equation*}
\text{ }\mathbf{H}(\mathbf{s;}\lambda ):=\left( \mathbf{F}(\mathbf{s;}%
\lambda );G(\mathbf{s;}\lambda )\right) =\left( F^{(1)}(\mathbf{s;}\lambda
),\ldots ,F^{(m)}(\mathbf{s;}\lambda );G(\mathbf{s;}\lambda )\right) ,
\end{equation*}%
and%
\begin{equation*}
\text{ }\mathbf{H}_{n}(\mathbf{s;}\lambda ):=\left( \mathbf{F}_{n}(\mathbf{s;%
}\lambda );G_{n+1}(\mathbf{s;}\lambda )\right) =\left( F_{n}^{(1)}(\mathbf{s;%
}\lambda ),\ldots ,F_{n}^{(m)}(\mathbf{s;}\lambda );G_{n+1}(\mathbf{s;}%
\lambda )\right),
\end{equation*}%
$n=0,1,2,...$ be a sequence of vector-valued m.p.g.f.'s selected from $\mathcal{F}%
^{m}\times \mathcal{G}$ in an iid manner in accordance with measure $\mathbb{%
Q}$. The sequence $\left\{ \mathbf{H}_{n},n\in \mathbb{N}_{0}\right\} $ is
called a \textit{random environment}. \ With each m.p.g.f. $F_{n}^{(i)}(%
\mathbf{s;}\lambda )$ we associate a random vector of offsprings $\mathbf{%
\xi }_{i}(n):=\left( \xi _{i1}(n),\xi _{i2}(n),\ldots ,\xi _{im}(n)\right) $
and a random variable $\varphi _{i}(n)$ such that%
\begin{equation*}
F_{n}^{(i)}(\mathbf{s;}\lambda ):=\mathbf{E}\left[ \mathbf{s}^{\mathbf{\xi }%
_{i}(n)}e^{-\lambda \varphi _{i}(n)}\right] =\mathbf{E}\left[ s_{1}^{\xi
_{i1}(n)}s_{2}^{\xi _{i2}(n)}\cdots s_{m}^{\xi _{im}(n)}e^{-\lambda \varphi
_{i}(n)}\right] \overset{d}{=}F^{(i)}(\mathbf{s;}\lambda )
\end{equation*}%
and with each m.p.g.f. $G_{n}(\mathbf{s;}\lambda )$ we associate a random
vector of immigrants $\mathbf{\eta }(n):=\left( \eta _{1}(n),\eta
_{2}(n),\ldots ,\eta _{m}(n)\right) $ and a random variable $\psi (n)$ such
that%
\begin{equation*}
G_{n}(\mathbf{s};\lambda ):=\mathbf{E}\left[ \mathbf{s}^{\mathbf{\eta }%
(n)}e^{-\lambda \psi (n)}\right] =\mathbf{E}\left[ s_{1}^{\eta
_{1}(n)}s_{2}^{\eta _{2}(n)}\cdots s_{m}^{\eta _{2}(n)}e^{-\lambda \psi (n)}%
\right] \overset{d}{=}G(\mathbf{s;}\lambda ).
\end{equation*}

Now we may give an informal description of  a multitype branching
processes with accumulation of a final product and immigration
which evolve in random environment (MBPFPIRE)
\begin{equation*}
\mathbf{T}(n)=\left( \mathbf{V}(n);\Theta (n)\right) ,\,n\in \mathbb{N}_{0},
\end{equation*}%
which may be treated as a process describing the evolution of a $m-$type
population of particles with immigration and accumulation of a final product.

Given an environment $\left\{ \mathbf{H}_{n},n\in \mathbb{N}_{0}\right\} $
the starting conditions for the process are: \ a vector of particles $%
\mathbf{V}(0)=\left( V_{1}(0),\ldots ,V_{m}(0)\right) $ (may be
random or equal zero), where $V_{i}(0)$ denotes the number of
particles of type $i\in \{1,\ldots ,m\}$ in the process at moment
$0$, and an amount $\Theta (0)$ (may be random or equal zero) of a
final product. All the initial particles have the unit life length
and just before the death produce children and some amount of the
final product independently of each other. A particle, say, of
type $i$, produces at the end of her life particles of different
types and adds some amount of the final product to the existing
amount of
the final product in accordance with m.p.g.f. $F_{0}^{(i)}(\mathbf{s;}%
\lambda )$. In addition, a random tuple of immigrants specified by the
vector $\mathbf{\eta }(1)=:\left( \eta _{1}(1),\eta _{2}(1),...,\eta
_{m}(1)\right) $ arrives to the system at moment $1$, where $\eta _{j}(1)$
is the number of type $j$ particles immigrating in the first generation of
the population, and the final product of size $\psi (1)$ is added to the
system with the joint distribution specified by the m.p.g.f. $G_{1}(\mathbf{s;}%
\lambda ).$ Thus, the amount of the final product $\Theta (1)$
accumulated in the system to this moment is
$$
\Theta (1)=\Theta (0)+\psi
(1)+\sum_{i=1}^m\sum_{k=1}^{V_i(0)}\varphi_i(0,k),
$$
where $\varphi_i(0,k)$ -- is the amount of the final product added
to the system at the death moment of the  $k$-particle of type $i$
among those which existed at moment $0.$

The newborn particles and immigrants entering the system at moment
$n=1$ constitute the first generation of the MBPIFPRE, have the
unit life-length and dying produce, independently of each other
and of the prehistory of the process, offsprings and
final product in accordance with their types and subject to the m.p.g.f. $%
F_{1}^{(i)}(\mathbf{s;}\lambda ),i=1,2,\ldots ,m.$ In addition, at
moment $n=2$ immigrants and some amount of final
product  described by a tuple $(%
\mathbf{\eta }(2);\psi (2))=:\left( \eta _{1}(2),\eta
_{2}(2),...,\eta
_{m}(2);\psi (2)\right)$ is contribute to the process, which is  specified by the m.p.g.f. $G_{2}(\mathbf{s;}%
\lambda ).$ And so on.

Note that in our model $\psi (n)$ may be positive even if $\mathbf{\eta }(n)=%
\mathbf{0}$.

\begin{definition}
\label{Definition1} A  $m$-type Galton-Watson branching process
\begin{equation*}
\mathbf{T}_{\varphi ,\psi }(n)=\mathbf{T}(n):=\left( \mathbf{V}(n);\Theta
(n)\right) =\left( V_{1}(n),...,V_{m}(n);\Theta (n)\right) ,\,n\in \mathbb{N}%
_{0}
\end{equation*}%
with immigration and  final product $\left( \varphi ,\psi \right)
$ in a fixed (but picked at random) environment $\left\{
\mathbf{H}_{n},\,n\in \mathbb{N}_{0}\right\} $ is a
time-inhomogeneous Markov process with the state space
\begin{equation*}
\mathbb{N}_{0}^{m}\times \mathbb{R}_{+}:=\left\{
z=(z_{1},...,z_{m};w),\,z_{i}\in \mathbb{N}_{0};w\in \lbrack 0,\infty
)\right\}
\end{equation*}%
defined as%
\begin{equation*}
\mathbf{T}(0)=\left( \mathbf{V}(0);\Theta (0)\right) =(\mathbf{v};\theta ),
\end{equation*}%
\begin{equation*}
\;\mathbf{E}\left[ \mathbf{s}^{\mathbf{V}(n+1)}e^{-\lambda \Theta (n+1)}|\;%
\mathbf{H}_{0},...,\mathbf{H}_{n};\mathbf{T}(0),,...,\mathbf{T}(n)\right]
=e^{-\lambda \Theta (n)}\left( \mathbf{F}_{n}(\mathbf{s;}\lambda )\right) ^{%
\mathbf{V}(n)}G_{n+1}(\mathbf{s;}\lambda ).
\end{equation*}
\end{definition}

Thus,%
\begin{eqnarray*}
\mathbf{T}(n+1) &=&\left( \mathbf{V}(n+1);\Theta (n+1)\right) \\
&=&\left( 0;\Theta (n)\right) +(\mathbf{\eta }(n+1);\psi
(n+1))+\sum_{i=1}^{m}\sum_{k=1}^{V_{i}(n)}\left( \mathbf{\xi }%
_{i}(n;k);\varphi _{i}(n;k)\right) ,
\end{eqnarray*}%
where $\left( \mathbf{\xi }_{i}(n;k);\varphi _{i}(n;k)\right) $ is
a random vector representing the offspring vector and the size of
the final product contributed to the process at the death moment
of the $k-$th particle of type $i$ of the $n-$th
generation. Given the environment and $\mathbf{V}(n)$ the tuple%
\begin{equation*}
\left( \mathbf{\xi }_{i}(n;k);\varphi _{i}(n;k)\right) ,\,k=1,2,\ldots
,V_{i}(n),i\in \{1,\ldots ,m\},\,(\mathbf{\eta }(n+1),\psi (n+1)),\text{ }%
n\in \mathbb{N}_{0};
\end{equation*}%
consists of independent vectors and, moreover, for each $n\in \mathbb{N}_{0}$
and $i\in \{1,\ldots ,m\}$ the vectors
\begin{equation*}
\left( \mathbf{\xi }_{i}(n;k);\varphi
_{i}(n;k)\right),\,k=1,2,\ldots ,V_{i}(n)
\end{equation*}%
are identically distributed: $\left( \mathbf{\xi }_{i}(n;k);\varphi
_{i}(n;k)\right) \overset{d}{=}\left( \mathbf{\xi }_{i}(n);\varphi
_{i}(n)\right) $.

In the sequel we write for brevity for fixed $n\in \mathbb{N}_{0}$%
\begin{equation*}
\mathbf{P}_{\mathbf{H}_{n}}\left( \mathbf{\cdot }\right) :=\mathbf{P}_{%
\mathbf{H}_{n}}\left( \mathbf{\cdot \,}|\,\mathbf{H}_{n}\right) ,\,\mathbf{E}%
_{\mathbf{H}_{n}}\left[ \mathbf{\cdot }\right] :=\mathbf{E}\left[ \mathbf{%
\cdot \,}|\,\mathbf{H}_{n}\right] ,
\end{equation*}%
\begin{equation*}
\mathbf{P}_{\mathbf{H}}\left( \mathbf{\cdot }\right) :=\mathbf{P}_{\mathbf{H}%
}\left( \mathbf{\cdot \,}|\,\mathbf{H}_{0},...,\mathbf{H}_{n}\right) \mathbf{%
,}\,\mathbf{E}_{\mathbf{H}}\left[ \mathbf{\cdot }\right] :=\mathbf{E}\left[
\mathbf{\cdot \,}|\,\mathbf{H}_{0},...,\mathbf{H}_{n}\right]
\end{equation*}%
or (where it will cause no confusion)
\begin{equation*}
\mathbf{P}_{\mathbf{H}}\left( \mathbf{\cdot }\right) :=\mathbf{P}_{\mathbf{H}%
}\left( \mathbf{\cdot \,}|\,\mathbf{H}_{0},...,\mathbf{H}_{n},...\right) ,\,%
\mathbf{E}_{\mathbf{H}}\left[ \mathbf{\cdot }\right] :=\mathbf{E}\left[
\mathbf{\cdot \,}|\,\mathbf{H}_{0},...,\mathbf{H}_{n},...\right]
\end{equation*}%
Similar meaning will have the notation $\mathbf{P}_{\mathbf{F}_{n}},\mathbf{E%
}_{\mathbf{F}_{n}},\mathbf{P}_{\mathbf{F}},\mathbf{E}_{\mathbf{F}}$ $\mathbf{%
\ }$and so on.

Observe that if $\Theta (0)=0$ and $\varphi _{i}(n;k)\equiv 1,\,
\psi (n)\equiv 0,\,n\in \mathbb{N}_{0},$ then $\Theta (N)-1$ is
the total
number of particles born and immigrating in the process within generations $%
0,1,...,N-1$; if%
\begin{equation*}
\varphi _{i}(n;k)=I\left\{ \sum_{j=1}^{m}\xi _{ij}(n;k)\geq
t\right\},\, \psi (n+1)\equiv 0,\,n\in \mathbb{N}_{0}
\end{equation*}%
for some positive integer $t$ (here and in what follows $I\left\{ A\right\} $
stands for the indicator of the event $A$) and $\Theta (0)=0,$ then $\Theta
(N)$ is the total number of particles of all types in generations $%
0,1,...,N-1$ each of which has at least $t$ children, and so on.

Letting $\lambda =0$ in Definition \ref{Definition1} we arrive to the
definition of the multitype branching process with immigration evolving in
random environment (MBPIRE) which we call the \textit{underlying} MBPIRE for
the initial MBPIFPRE.

\begin{definition}
A $m$-type Galton-Watson branching process with immigration
$\mathbf{V}(n)$
in a fixed (but picked at random) environment $\left\{ (\mathbf{f}%
_{n},g_{n+1}),n\in \mathbb{N}_{0}\right\} $ is a time-inhomogeneous Markov
chain with the state space
\begin{equation*}
\mathbb{N}_{0}^{m}:=\left\{ \mathbf{z}=(z_{1},...,z_{m}),z_{i}\in \mathbb{N}%
_{0}\right\}
\end{equation*}%
defined as $\mathbf{V}(0)=\mathbf{z}$ and for $n\in \mathbb{N}_{0}$
\begin{equation}
\;\mathbf{E}\left[ \mathbf{s}^{\mathbf{V}(n+1)}|\;\mathbf{f}_{0},g_{1},...,%
\mathbf{f}_{n},g_{n+1};\mathbf{V}(0),...,\mathbf{V}(n)\right] =\left(
\mathbf{f}_{n}\left( \mathbf{s}\right) \right) ^{\mathbf{V}(n)}g_{n+1}(%
\mathbf{s}).  \label{Immig1}
\end{equation}
\end{definition}

Thus,%
\begin{equation*}
\mathbf{V}(n+1)=\sum_{i=1}^{m}\sum_{k=1}^{V_{i}(n)}\mathbf{\xi }_{i}(n;k)+%
\mathbf{\eta }(n+1),
\end{equation*}%
where, given the environment and fixed $\mathbf{V}(n),$ the tuple
\begin{equation*}
\mathbf{\xi }_{i}(n;k),k=1,2,...V_{i}(n);i=1,2,...,m;\,\mathbf{\eta }(n+1)
\end{equation*}%
consists of independent random variables.

In the sequel we need the definition of the ordinary multitype branching
process with final product which evolves in random environment (MBPFPRE).

\begin{definition}
A $m$-type Galton-Watson branching process
\begin{equation*}
\mathbf{R}_{\varphi }(n)=\mathbf{R}(n):=\left( \mathbf{Z}(n);\Phi (n)\right)
=\left( Z_{1}(n),\ldots ,Z_{m}(n);\Phi (n)\right) ,\,n\in \mathbb{N}_{0},
\end{equation*}%
with final product $\varphi $ in a fixed (but picked at random) environment $%
\left\{ \mathbf{F}_{n},\,n\in \mathbb{N}_{0}\right\} $ is a
time-inhomogeneous Markov process with the state space
\begin{equation*}
\mathbb{N}_{0}^{m}\times \mathbb{R}_{0}:=\left\{ z=(z_{1},\ldots
,z_{m};w),\,z_{i}\in \mathbb{N}_{0};w\in \lbrack 0,\infty )\right\}
\end{equation*}%
defined as%
\begin{equation*}
\mathbf{R}(0)=\left( \mathbf{Z}(0);\Phi (0)\right) =(\mathbf{z};\varphi _{0})%
\mathbf{,}
\end{equation*}%
\begin{equation*}
\;\mathbf{E}\left[ \mathbf{s}^{\mathbf{Z}(n+1)}e^{-\lambda \Phi (n+1)}|\;%
\mathbf{F}_{0},\ldots ,\mathbf{F}_{n};\mathbf{R}(0),\ldots ,\mathbf{R}(n)%
\right] =e^{-\lambda \Phi (n)}\left( \mathbf{F}_{n}(\mathbf{s;}\lambda
)\right) ^{\mathbf{Z}(n)}.
\end{equation*}
\end{definition}

Note that the initial value $\mathbf{R}(0)$ may be random. Thus,
\begin{equation}
\Phi (n+1)=\Phi (0)+\sum_{l=0}^{n}\sum_{i=1}^{m}\sum_{k=1}^{Z_{i}(l)}\varphi
_{i}(l;k)  \label{DDphi}
\end{equation}%
where $\left( \mathbf{\xi }_{i}(n;k);\varphi _{i}(n;k)\right) $ is a random
vector representing the offspring vector and the size of the final product
of the $k-$th particle of type $i$ of the $n-$th generation of the process.

Finally, excluding immigration and final product we obtain the
definition of the ordinary \textit{underlying }multitype branching
process in random environment.

\begin{definition}
A $m$-type Galton-Watson process
\begin{equation*}
\mathbf{Z}(n)=\left( Z_{1}(n),\ldots ,Z_{m}(n)\right) ,\,n\in \mathbb{N}_{0}
\end{equation*}%
in a fixed (but selected at random) environment $\left\{ \mathbf{f}%
_{n},\,n\in \mathbb{N}_{0}\right\} $ is a time-inhomogeneous Markov chain
with the state space
\begin{equation*}
\mathbb{N}_{0}^{m}:=\left\{ z=(z_{1},\ldots ,z_{m}),\,z_{i}\in \mathbb{N}%
_{0}\right\}
\end{equation*}%
defined as
\begin{equation}
\mathbf{Z}(0)=\mathbf{z},\;\mathbf{E}\left[ \mathbf{s}^{\mathbf{Z}(n+1)}|\;%
\mathbf{f}_{0},\ldots ,\mathbf{f}_{n};\mathbf{Z}(0),\ldots ,\mathbf{Z}(n)%
\right] =\left( \mathbf{f}_{n}\left( \mathbf{s}\right) \right) ^{\mathbf{Z}%
(n)}.  \label{defVP}
\end{equation}
\end{definition}

It is necessary to mention that the one-dimensional BPRE without
immigration was introduced by Smith and Wilkinson in \cite{SW68}
and, in a more general setting in \cite{AK71a} and \cite{AK71} and
has been investigated by many authors (see survey \cite{VZ} for a
list of references up to 1985 and
\cite{AGKV},\cite{AGKV2},\cite{Koz},
\cite{VD}-\cite{VD1},\cite{GDV},\cite{GK00}, and \cite{KOZ95} for
some more recent results).

The ordinary single-type Galton-Watson branching processes with final
product where investigated by Sevastyanov \cite{Se72} \ for the particular
case $\varphi (n;k)\equiv 1,$ and by Grishechkin \cite{Gri90} (for the
general $\varphi (n;k)$). MBPRE were analyzed in \cite{AK71a} and \cite%
{Tan81} and MBPIRE were considered by \cite{Key87} and \cite{Roi07}. MBPFPRE
were introduced in \cite{Vat2009} to study the tail distribution of busy
periods of the branching type polling systems with zero switch-over times.
In the present paper we analyze MBPIFPRE and apply the obtained results to
investigate properties of busy periods of the branching type polling systems
with random and, generally speaking, positive switch-over times.

Let $\mathbf{e}_{i}:=\left( 0,\ldots ,1,\ldots ,0\right) $ be a $m-$%
dimensional vector with zero components except the $i$-th equal to $1,$ $%
\mathbf{0}$ and $\mathbf{1}$ are $m-$dimensional vectors all whose
components are equal to $0$ and $1$, respectively.

Denote
\begin{equation}
A_{n}=\left( a_{ij}(n)\right) _{i,j=1}^{m}:=\left( \frac{\partial
f_{n}^{(i)}(\mathbf{s})}{\partial s_{j}}\Big| _{\mathbf{s}=\mathbf{1}%
}\right) _{i,j=1}^{m}  \label{meanMat}
\end{equation}%
the mean matrix of the vector-valued p.g.f. $\mathbf{f}_{n}$ ,%
\begin{equation}
\mathbf{C}_{n}:=\left( \mathbf{E}_{\mathbf{F}_{n}}\varphi _{1}(n),...,%
\mathbf{E}_{\mathbf{F}_{n}}\varphi _{1}(n)\right) ^{\prime }=\left( \frac{%
dF_{n}^{(1)}(\mathbf{1},\lambda )}{d\lambda }\Big|_{\lambda
=0},...,\frac{dF_{n}^{(m)}(\mathbf{1},\lambda )}{d\lambda }\,\Big|
_{\lambda =0}\right) ^{\prime }  \label{finprod}
\end{equation}%
the mean vector of the amount of the final product of particles of the $n$%
-th generation and%
\begin{equation}
\mathbf{B}_{n}:=\left( B_{1}(n),...,B_{m}(n)\right) ^{\prime }:=\left( \frac{%
\partial g_{n}(\mathbf{s})}{\partial s_{1}}\Big| _{\mathbf{s}=\mathbf{1}%
} ,...,\frac{\partial g_{n}(\mathbf{s})}{\partial s_{m}}\Big| _{%
\mathbf{s}=\mathbf{1}}\right) ,  \label{meanVector}
\end{equation}%
\begin{equation}
D_{n}:=\mathbf{E}_{G_{n}}\psi (n)=\frac{dG_{n}(\mathbf{1},\lambda )}{%
d\lambda }\Big|_{\lambda =0} \label{FinImproduct}
\end{equation}%
the respective characteristics for immigrants.

By our assumptions the tuples $\left( A_{n},\mathbf{B}_{n+1},\mathbf{C}%
_{n},D_{n+1}\right) ,\,n\in \mathbb{N}_{0}$ are iid:
\begin{equation*}
\left( A_{n},\mathbf{B}_{n+1},\mathbf{C}_{n},D_{n+1}\right) \overset{d}{=}(A,%
\mathbf{B},\mathbf{C},D).
\end{equation*}

For vectors $\mathbf{u}=(u_{1},\ldots ,u_{m})$ and $\mathbf{v=}(v_{1},\ldots
,v_{m})^{\prime }\in \mathbb{R}^{m}$ denote
\begin{equation*}
\left\langle \mathbf{u},\mathbf{v}\right\rangle :=\sum_{k=1}^{m}u_{i}v_{i}
\end{equation*}%
their inner product.

For a $m\times m$ \ matrix $A=\left( a_{ij}\right) _{i,j=1}^{m}$ and a $m-$%
dimensional vector $\mathbf{u=}(u_{1},\ldots ,u_{m})$ introduce the norms
\begin{equation*}
\left\Vert A\right\Vert :=\sum_{i,j=1}^{m}\left\vert a_{ij}\right\vert
,\quad \left\Vert \mathbf{u}\right\Vert :=\sum_{i=1}^{m}\left\vert
u_{i}\right\vert
\end{equation*}%
and%
\begin{equation*}
\left\Vert A\right\Vert _{2}:=\sqrt{\sum_{i,j=1}^{m}\left\vert
a_{ij}\right\vert ^{2}},\quad \left\Vert \mathbf{u}\right\Vert _{2}:=\sqrt{%
\sum_{i=1}^{m}\left\vert u_{i}\right\vert ^{2}}.
\end{equation*}

Let, further,
\begin{equation*}
\Pi _{l,n}:=\prod_{i=l}^{n-1}A_{i},\qquad 1\leq l\leq n,
\end{equation*}%
with the agreement that $\Pi _{n,n}:=E$ is the unit $m\times m$ matrix.

\section{\protect\bigskip Statement of main results}

Denote $\mathcal{D}:=\left\{ x>0:\mathbf{E}\left\Vert A\right\Vert
^{x}<\infty \right\} $ and, for a given $x\geq 0$, set%
\begin{equation}
s(x):=\lim_{n\rightarrow \infty }\left( \mathbf{E}\left\Vert A_{n-1}\cdots
A_{0}\right\Vert ^{x}\right) ^{1/n}=\lim_{n\rightarrow \infty }\left(
\mathbf{E}\left\Vert \Pi _{0,n}\right\Vert ^{x}\right) ^{1/n}  \label{defK}
\end{equation}%
and let%
\begin{equation}
\alpha :=s^{\prime }(0)=\lim_{n\rightarrow \infty }\frac{1}{n}\mathbf{E}\log
\left\Vert A_{n-1}\cdots A_{0}\right\Vert =\lim_{n\rightarrow \infty }\frac{1%
}{n}\mathbf{E}\log \left\Vert \Pi _{0,n}\right\Vert  \label{DefK0}
\end{equation}%
be the top Lyapunov exponent for this sequence of matrices.

It is known that the limits in (\ref{defK}) and (\ref{DefK0})
exist and, moreover, $s(x)$ is a log-convex continuous function in
$\mathcal{D}$ (see, for instance, \cite{Kin73}). Put
\begin{equation}
\kappa :=\inf \left\{ x>0:s(x)>1\right\}  \label{DefKappa}
\end{equation}%
and $\kappa =\infty $ if $s(x)\leq 1$ for all $x>0$. Observe that $s(0)=1$
and, therefore, $\kappa =0$ if $\alpha >0$ and $\kappa \in (0,\infty ]$ if $%
\alpha <0$.

In the last case (which will be our main concern) the series $%
\sum_{n=0}^{\infty }\mathbf{E}\left\Vert \Pi _{0,n}\right\Vert ^{x}$
converges if $0<x<\kappa $ and diverges if $x>\kappa $.

In what follows we call the \textit{underlying } MBPRE subcritical if $%
\alpha <0$ and supercritical if $\alpha >0$. The same terminology we keep
for the respective MBPIRE.

Now we formulate important statements concerning properties of MBPRE and
MBPIRE.

Let%
\begin{equation*}
q_{i}\left( \mathbf{f}\right) :=\lim_{n\rightarrow \infty }\mathbf{P}_{%
\mathbf{f}}\left( \left\Vert \mathbf{Z}(n)\right\Vert =0|\mathbf{Z}(0)=%
\mathbf{e}_{i}\right)
\end{equation*}%
be the extinction probability of a MBPRE initiated at time $0$ by a single
individual of type $i$ and
\begin{equation*}
\mathbf{q}\left( \mathbf{f}\right) :=\left( q_{1}\left( \mathbf{f}\right)
,\ldots,q_{m}\left( \mathbf{f}\right) \right) .
\end{equation*}

\begin{theorem}
\label{Ttan}(\cite{Tan81}) If the mean matrices of a MBPRE meets the
condition
\begin{equation}
\mathbf{E}\log ^{+}\left\Vert A\right\Vert <\infty  \label{ExpCond}
\end{equation}%
and there exists a positive integer $L$ such that%
\begin{equation*}
\mathbf{P}\left( \min_{1\leq i,j\leq m}\left( A_{L-1}A_{L-2}\cdots
A_{0}\right) _{ij}>0\right) =1
\end{equation*}%
and $1\leq l\leq m$ such that%
\begin{equation*}
\mathbf{E}\left\vert \log \left( 1-\mathbf{P}_{\mathbf{f}}\left( \mathbf{Z}%
_{l}(L)=0\,|\,\mathbf{Z}(0)=\mathbf{e}_{i}\right) \right) \right\vert
<\infty ,
\end{equation*}%
then, for $\alpha $ specified by (\ref{DefK0})

1) $\alpha <0$ implies $\mathbf{P}\left( \mathbf{q}\left( \mathbf{f}\right) =%
\mathbf{1}\right) =1;$

2) $\alpha >0$ implies $\mathbf{P}\left( \mathbf{q}\left( \mathbf{f}\right) <%
\mathbf{1}\right) =1$ and%
\begin{equation}
\mathbf{P}_{\mathbf{f}}\left( \lim_{n\rightarrow \infty }n^{-1}\log
\left\Vert \mathbf{Z}(n)\right\Vert =\alpha \,|\,\mathbf{Z}(0)=\mathbf{e}%
_{i}\right) =1-q_{i}\left( \mathbf{f}\right)  \label{SU}
\end{equation}%
with probability 1 for $1\leq i\leq m.$
\end{theorem}

The next statement deals with MBPIRE.

\begin{theorem}
\label{Key}(\cite{Key87} and \cite{Roi07}) Let a MBPIRE satisfy condition (%
\ref{ExpCond}), $\alpha <0$ and
\begin{equation*}
\mathbf{E}\log ^{+}\left\Vert \mathbf{B}\right\Vert <\infty .
\end{equation*}%
Then, for any $\mathbf{v}\in \mathbb{N}_{0}^{m}$ the limit%
\begin{equation*}
\lim_{n\rightarrow \infty }\mathbf{P}\left( \mathbf{V}(n)=\mathbf{v}\right)
=:D(\mathbf{v})
\end{equation*}%
exists and defines a probability distribution on $\mathbb{N}_{0}^{m}$. \ If,
in addition, there exists $L\geq 1$ such that
\begin{equation*}
\mathbf{P}\left( \mathbf{P}_{\mathbf{g}}\left( \mathbf{V}(L)=\mathbf{0|V}(0)=%
\mathbf{1}\right) >0\right) >0
\end{equation*}%
then $D(\mathbf{0})>0$ and, therefore,
\begin{equation*}
\mathbf{P}\left( \mathbf{V}(n)=\mathbf{0}\,\text{ i.o. }\right)
=1.
\end{equation*}
\end{theorem}

Introduce the set
\begin{equation*}
U_{+}=\left\{ \mathbf{u}=(u_{1},\ldots ,u_{m})\in \mathbb{R}^{m}:u_{i}\geq
0,1\leq i\leq m,\,\left\Vert \mathbf{u}\right\Vert _{2}=1\right\}
\end{equation*}%
and associate with the tuple $(A_{n},\mathbf{C}_{n}\,),$ $n\in \mathbb{N}%
_{0} $ of iid pairs the series of vectors
\begin{equation}
\Xi _{l}:=\sum_{k=l}^{\infty }A_{l}A_{l+1}\cdots A_{k-1}\mathbf{C}%
_{k}=\sum_{k=l}^{\infty }\Pi _{l,k}\mathbf{C}_{k},\ l\,\in \mathbb{N}%
_{0};\,\quad \Xi :=\Xi _{0}.  \label{DefSumThet}
\end{equation}

Our main results are established under the following hypothesis.

\textbf{Condition T.} There exist positive constants $\kappa $ and $K_{0}$
and a continuous strictly positive function $l(\mathbf{u})$ on $U_{+}$ such
that for all $\mathbf{u\in }U_{+}$
\begin{equation}
\lim_{y\rightarrow \infty }y^{\kappa }\mathbf{P}\left( \left\langle \mathbf{%
u,}\Xi \right\rangle >y\right) =K_{0}l(\mathbf{u}).  \label{FcondT}
\end{equation}

In Section \ref{SecAux} we list sufficient conditions imposed on the
distributions of the pairs $(A_{n},\mathbf{C}_{n})$ which provide the
validity of Condition $T$. These conditions are extracted from paper \cite%
{Kes74} where the behavior of the tail distribution of sums and products of
random matrices were investigated.

Let
\begin{equation}
\Phi :=\lim_{n\rightarrow \infty }\Phi (n)=\Phi (0)+\sum_{n=0}^{\infty
}\sum_{i=1}^{m}\sum_{k=1}^{Z_{i}(n)}\varphi _{i}(n;k)  \label{DefSums}
\end{equation}%
be the total size of the final product produced by the particles of the
MBPFPRE up to the extinction moment (if any).

The following theorem was proved in \cite{Vat2009}.

\begin{theorem}
\label{TailMain}Let a MBPFPRE satisfy the following hypotheses:

1) the underlying MBPRE \ is subcritical and meets the conditions of Theorem~%
\ref{Ttan};

2) for $\kappa $ specified by (\ref{DefKappa}) the following assumptions
fulfill:

if $\kappa >1$ then
\begin{equation}
\max_{1\leq i\leq m}\mathbf{E}\left\vert \sum_{j=1}^{m}\left( \xi _{ij}-%
\mathbf{E}_{\mathbf{f}}\xi _{ij}\right) \right\vert ^{\kappa }<\infty \text{
and \ \ }\mathbf{E}\left\vert \sum_{i=1}^{m}\left( \varphi _{i}-\mathbf{E}_{%
\mathbf{F}}\varphi _{i}\right) \right\vert ^{\kappa }<\infty ,
\label{Ckappabig}
\end{equation}

if $\kappa \leq 1$ then%
\begin{equation}
\max_{1\leq i\leq m}\mathbf{E}\left( \sum_{j=1}^{m}Var_{\mathbf{f}}\xi
_{ij}\right) ^{\kappa }<\infty \text{ and \ }\mathbf{E}\left(
\sum_{i=1}^{m}Var_{\mathbf{F}}\varphi _{i}\right) ^{\kappa }<\infty ;
\label{Ckappasmall}
\end{equation}

4) there exists $\delta >0$ such that $0<\mathbf{E}\varphi _{i}^{\kappa
+\delta }<\infty ,\quad i=1,\ldots ,m;$

5) the mean matrix (\ref{meanMat}) and the vector (\ref{finprod}) are such
that Condition $T$ is valid.

Then, as $y\rightarrow \infty $%
\begin{equation}
\mathbf{P}\left( \Phi >y\,|\,\mathbf{Z}(0)=\mathbf{z}\right) \sim C(\mathbf{z%
})y^{-\kappa },\qquad C(\mathbf{z})\in (0,\infty ).  \label{SimplFin}
\end{equation}
\end{theorem}

The theorem has been used in \cite{Vat2009} to investigate the
tail distributions of the busy periods and some other
characteristics for a wide class of polling systems with zero
switch-over time. In the present paper we consider MBPIFPRE and
apply in Section \ref{Spoll} the respective results to deduce
similar statements for busy periods and certain other
characteristics of polling systems with positive (possibly random)
switch-over times. To this aim we introduce the notion of life
periods of a MBPIFPRE.

\begin{definition}
We say that a branching process $\mathbf{T}(n),n\in
\mathbb{N}_{0}$ with $m$ types of particles, immigration and final
product has a life period of
length $\Upsilon _{N}$ initiated at moment $N\geq 1$ if $\left( \mathbf{V}%
(N-1);\Theta (N-1)\right) =(0;\theta ),$ $\left\Vert \mathbf{\eta }%
(N)\right\Vert >0,$ and%
\begin{equation*}
\Upsilon _{N}:=\min \left\{ k>N:\mathbf{V}(k)=0\right\} .
\end{equation*}
\end{definition}

In what follows we consider without loss of generality the
distributions of various characteristics related with the life
period of length $\Upsilon
:=\Upsilon _{1}$, and study the tail distribution of the quantity%
\begin{equation*}
\Theta =\Theta (\Upsilon ):=\theta +\sum_{n=0}^{\Upsilon -1}\left( \psi
(n+1)+\sum_{i=1}^{m}\sum_{k=1}^{V_{i}(n)}\varphi _{i}(n;k)\right)
\end{equation*}%
- the total amount of the final product accumulated in the MBPIFPRE during
the life-period which starts by $\mathbf{\eta }(1)=\mathbf{z}\neq \mathbf{0}$
particles at moment $N=1$. Since our proofs (but not constants!) do not
depend on the particular value of $\mathbf{\eta }(1)$ we do not specify it
explicitly in the subsequent arguments.

Note that if $\theta =0$ and $\varphi _{i}(n;k)\equiv 1,$ $\psi (n)\equiv 0$
for all $n\in \mathbb{N}_{0}$ then $\Theta $ is the total number of
individuals existing in the MBPIFPRE during the respective life-period.

Now we are ready to formulate the main result of the present
paper.

\begin{theorem}
\label{TTmain}Let the conditions of Theorem \ref{TailMain} be
valid,
\begin{equation*}
\max_{1\leq i\leq m}\mathbf{E}\eta _{i}^{1+\kappa }<\infty ,\qquad \mathbf{E}%
\left\vert \psi -\mathbf{E}_{\mathbf{G}}\psi \right\vert ^{\kappa
}<\infty
\end{equation*}%
and, in addition,
\begin{equation*}
\mathbf{E}\left\Vert \mathbf{B}\right\Vert ^{q}<\infty \text{ \ \ for some }%
q>0\text{ and \ \ \ }\mathbf{E}\left\Vert A\right\Vert ^{\kappa }<\infty .
\end{equation*}%
Then, as $y\rightarrow \infty $%
\begin{equation}
\mathbf{P}\left( \Theta >y\right) \sim C_{I}y^{-\kappa },\qquad C_{I}\in
(0,\infty ).  \label{EstTail}
\end{equation}%
In particular,%
\begin{equation}
\mathbf{E}\Theta ^{x}<\infty  \label{EstMom}
\end{equation}%
if and only if $x<\kappa .$
\end{theorem}

In  Section \ref{Spoll} we use Theorem  \ref{TTmain}  to show that
 the tail distributions of the busy periods of a
wide class of branching type polling systems whose input
parameters and service disciplines vary in a random manner and
switch-over times  are positive, decay at infinity as
$y^{-\kappa}$ for some $\kappa>0.$

\section{Auxiliary results\label{SecAux}}

The proof of Theorem \ref{TTmain} uses Condition $T$ whose validity is not
easy to check. We list here a set of assumptions given in \cite{Kes74} which
imply Condition T.

Let $\Lambda (A)$ be the spectral radius of the matrix $A$. The following
statement is a refinement of a Kesten theorem from \cite{Kes74}.

\begin{theorem}
\label{Tkest}(see \cite{DMB99}) Let $\left\{ A_{n},n\geq 0\right\} $ be a
sequence of iid matrices generated by a measure $\mathbb{P}_{A}$ with
support concentrated on nonnegative matrices and $A=\left( a_{ij}\right)
_{i,j=1}^{m}\overset{d}{=}A_{n}.$ Assume that the following conditions are
valid:

1) there exists $q>0$ such that $E\left\Vert A\right\Vert ^{q}<\infty $;

2) $A$ has no zero rows a.s.;

3) the group generated by
\begin{equation*}
\left\{ \log \Lambda (a_{n}\cdots a_{0}):a_{n}\cdots a_{0}>0\text{ for some }%
n\text{ and }a_{i}\in \text{supp}(\mathbb{P}_{A})\right\}
\end{equation*}%
is dense in $\mathbb{R};$

4) there exists $\kappa _{0}>0$ for which
\begin{equation*}
\mathbf{E}\left[ \min_{1\leq i\leq m}\left( \sum_{j=1}^{m}a_{ij}\right)
^{\kappa _{0}}\right] \geq m^{\kappa _{0}/2}
\end{equation*}%
and%
\begin{equation*}
\mathbf{E}\left\Vert A\right\Vert ^{\kappa _{0}}\log ^{+}\left\Vert
A\right\Vert <\infty .
\end{equation*}

Then there exists a $\kappa \in (0,\kappa _{0}]$ such that%
\begin{equation*}
s^{\prime }(\kappa )=\lim_{n\rightarrow \infty }\frac{1}{n}\log \mathbf{E}%
\left\Vert A_{n-1}\cdots A_{0}\right\Vert ^{\kappa }=0.
\end{equation*}

If, in addition, the tuple of $m$-dimensional vectors $\left\{ \mathbf{C}%
_{n},n\geq 0\right\} $ is such that the pairs $(A_{n},\mathbf{C}_{n})$, $%
n\in \mathbb{N}_{0}$ are iid: $\left( A_{n},\mathbf{C}_{n}\right) \overset{d}%
{=}\left( A,\mathbf{C}\right) $ and such that
\begin{equation*}
\mathbf{P}\left( \mathbf{C}=\mathbf{0}\right) <1,\quad \mathbf{P}\left(
\mathbf{C}\geq \mathbf{0}\right) =1,\quad \mathbf{E}\left\Vert \mathbf{C}%
\right\Vert ^{\kappa }<\infty ,
\end{equation*}%
then there exist a constant $K_{0}\in \left( 0,\infty \right) $ and a
continuous strictly positive function $l(\mathbf{u})$ on $U_{+}$ such that
relation (\ref{FcondT}) holds.
\end{theorem}

The next lemma provides estimates of moments of a random walk (see, for
instance, \cite{Gut}, Theorem 1.5.1):

\begin{lemma}
\label{LGut}If $X_{i},i=1,2,...$ is a sequence of iid random variables such
that $\mathbf{E}\left\vert X_{i}\right\vert ^{p}<\infty $ and $\mathbf{E}%
X_{i}=0$ if $p\geq 1$ and $M$ is a stopping time for the sequence $\Gamma
_{n}:=X_{1}+...+X_{n}$ then there exists a constant $R_{p}\in (0,\infty )$
such that
\begin{equation*}
\mathbf{E}\left\vert \Gamma _{M}\right\vert ^{p}\leq R_{p}\mathbf{E}%
\left\vert X_{i}\right\vert ^{p}\mathbf{E}M^{p/2\vee 1}.
\end{equation*}
\end{lemma}

\bigskip From now on we agree to denote by $K,K_{1},K_{2},...$ positive
constants which may be different from formula to formula.

\begin{lemma}
\label{LGeom}(\cite{Vat2009}) If a subcritical MBPRE starts by $\mathbf{Z}%
(0)=\mathbf{z}$ individuals and parameter $\kappa $ in (\ref{DefKappa})
exceeds $1$ then for each $x\in \lbrack 1,\kappa )$ there exist $\rho =\rho
(x)\in \left( 0,1\right) $ and $K=K(x)<\infty $ such that%
\begin{equation}
\mathbf{E}\left\Vert \mathbf{Z}(n)\right\Vert ^{x}\leq K\rho ^{n}\left\Vert
\mathbf{z}\right\Vert ^{x}  \label{geom11}
\end{equation}%
for all $n\in \mathbb{N}_{0}$.
\end{lemma}

The next statement is borrowed from \cite{Key87}, page 351.

\begin{lemma}
\label{LtailAsympt}If conditions of Theorem \ref{Key} are valid and, in
addition,
\begin{equation*}
\mathbf{E}\left\Vert \mathbf{B}\right\Vert ^{q}<\infty \text{ \ and \ \ \ }%
\mathbf{E}\left\Vert A\right\Vert ^{q}<\infty \text{ \ \ for some }q>0,
\end{equation*}%
then there exist positive constants $K_{1}$ and $K_{2}$ such that for all $%
t>0$%
\begin{equation*}
\mathbf{P}\left( \Upsilon >t\right) \leq K_{1}e^{-K_{2}t}.
\end{equation*}
\end{lemma}

Let
\begin{equation*}
\mathbf{Z}^{I}(i,j;n,k):=\left(
Z_{1}^{I}(i,j;n,k),...,Z_{m}^{I}(i,j;n,k)\right)
\end{equation*}%
be the $m-$dimensional vector of the total progeny alive at moment
$k>n$ of the $j-$th immigrant of type $i$ entering the process at
time $n$,
\begin{equation*}
\mathbf{Z}^{I}(n,k):=\sum_{i=1}^{m}\sum_{j=1}^{\eta _{i}(n)}\mathbf{Z}%
^{I}(i,j;n,k)
\end{equation*}%
be the $m-$dimensional vector of the total progeny alive at moment
$k>n$ of all the immigrants entering the process at time $n$, and
\begin{equation*}
\mathbf{Y}(n):=\mathbf{\eta }(n)+\sum_{k=n+1}^{\infty }\mathbf{Z}^{I}(n,k)
\end{equation*}%
be the $m-$dimensional vector representing the total progeny
steaming from the immigrants entering the process at time $n$.

\begin{lemma}
\label{LEstzero}If the conditions of Theorem \ref{TailMain} are valid and,
in addition,
\begin{equation*}
\max_{1\leq i\leq m}\mathbf{E}\eta _{i}^{\kappa +1}<\infty
\end{equation*}%
then there exists a constant $K_{1}\in \left( 0,\infty \right) $ such that
for all $x>0$
\begin{equation*}
\mathbf{P}\left( \left\Vert \mathbf{Y}(1)\right\Vert \geq x\right) \leq
K_{1}x^{-\kappa }.
\end{equation*}
\end{lemma}

\textbf{Proof}. We have
\begin{equation*}
\mathbf{Y}(1)=\mathbf{\eta }(1)+\sum_{k=1}^{\infty
}\sum_{i=1}^{m}\sum_{j=1}^{\eta _{i}(1)}\mathbf{Z}^{I}(i,j;1,k+1).
\end{equation*}%
Hence, observing that, for any fixed $i\in \left\{ 1,...,m\right\} $ and any
$j=1,2,...,\eta _{i}(1)$
\begin{equation*}
\left\{ \mathbf{Z}^{I}(i,j;1,k+1),k=1,2,...\right\}
\overset{d}{=}\left\{
\mathbf{Z}(k),k=1,2,...\,|\,\mathbf{Z}(0)=\mathbf{e}_{i}\right\}
\end{equation*}%
and using Theorem \ref{TailMain} with $\varphi _{i}\equiv 1,\,i\in \left\{
1,...,m\right\} $ we see that
\begin{eqnarray*}
\mathbf{P}\left( \left\Vert \mathbf{Y}(1)\right\Vert \geq x\right) &\leq &%
\mathbf{P}\left( \left\Vert \mathbf{\eta }\right\Vert \geq \frac{x}{2}%
\right) +\sum_{i=1}^{m}\mathbf{P}\left( \sum_{j=1}^{\eta _{i}(1)}\left\Vert
\sum_{k=1}^{\infty }\mathbf{Z}^{I}(i,j;1,k+1)\right\Vert \geq \frac{x}{2m}%
\right) \\
&\leq &\frac{2^{\kappa }\mathbf{E}\left\Vert \mathbf{\eta }\right\Vert
^{\kappa }}{x^{\kappa }}+\sum_{i=1}^{m}\mathbf{E}\left[ \sum_{j=1}^{\eta
_{i}(1)}\mathbf{P}\left( \sum_{k=1}^{\infty }\left\Vert \mathbf{Z}%
(k)\right\Vert \geq \frac{x}{m\eta _{i}(1)}\big|\,\mathbf{Z}(0)=\mathbf{e}%
_{i};\eta _{i}(1)\right) \right] \\
&=&\frac{2^{\kappa }\mathbf{E}\left\Vert \mathbf{\eta }\right\Vert ^{\kappa }%
}{x^{\kappa }}+\sum_{i=1}^{m}\mathbf{E}\left[ \eta _{i}(1)\mathbf{P}\left(
\sum_{k=1}^{\infty }\left\Vert \mathbf{Z}(k)\right\Vert \geq \frac{x}{m\eta
_{i}(1)}\big|\,\mathbf{Z}(0)=\mathbf{e}_{i};\eta _{i}(1)\right) \right] \\
&\leq &\frac{2^{\kappa }\mathbf{E}\left\Vert \mathbf{\eta }\right\Vert
^{\kappa }}{x^{\kappa }}+K_{1}\sum_{i=1}^{m}\mathbf{E}\left[ \eta _{i}(1)%
\frac{\left( m\eta _{i}(1)\right) ^{\kappa }}{x^{\kappa }}\right] \leq \frac{%
K_{2}}{x^{\kappa }}\sum_{i=1}^{m}\mathbf{E}\eta _{i}^{\kappa +1}=\frac{K_{1}%
}{x^{\kappa }}.
\end{eqnarray*}

The lemma is proved.

Let%
\begin{equation*}
\mu =\mu (r)=\min \left\{ n:\left\Vert \mathbf{V}(n)\right\Vert >r\right\}
\end{equation*}%
be the first moment when the size of the population of the MBPIRE
exceeds level $r$ with the natural agreement that $\mu =\infty $
if $\max_{n}\left\Vert \mathbf{V}(n)\right\Vert \leq r$.

\begin{lemma}
\label{Lke7}Under the conditions of Theorem \ref{TTmain} for any $%
\varepsilon >0$ there exists $r_{1}=r_{1}(\varepsilon )>0$ such that%
\begin{equation*}
\mathbf{P}\left( \sum_{n=\mu +1}^{\Upsilon -1}\left\Vert \mathbf{Y}%
(n)\right\Vert \geq \varepsilon y,\mu <\Upsilon \right) \leq \varepsilon
y^{-\kappa }
\end{equation*}%
for all $r\geq r_{1}.$
\end{lemma}

\textbf{Proof}. In view of the estimate
\begin{equation}
\sum_{n=1}^{\infty }n^{-2}=\pi ^{2}/6\leq 2  \label{Series}
\end{equation}%
we have
\begin{eqnarray}
&&\mathbf{P}\left( \sum_{n=\mu +1}^{\Upsilon -1}\left\Vert \mathbf{Y}%
(n)\right\Vert \geq \varepsilon y,\mu <\Upsilon \right)  \notag \\
&&\qquad=\mathbf{P}\left( \sum_{n=1}^{\infty }\left\Vert
\mathbf{Y}(n)\right\Vert I\left\{ \mu <n<\Upsilon \right\} \geq
6\pi ^{-2}\varepsilon
y\sum_{n=1}^{\infty }\frac{1}{n^{2}}\right)  \notag \\
&&\qquad\leq \sum_{n=1}^{\infty }\mathbf{P}\left( \left\Vert \mathbf{Y}%
(n)\right\Vert I\left\{ \mu <n<\Upsilon \right\} \geq \frac{\varepsilon y}{%
2n^{2}}\right)  \notag \\
&&\qquad=\sum_{n=1}^{\infty }\mathbf{P}\left( \mu <n<\Upsilon \right) \mathbf{P}%
\left( \left\Vert \mathbf{Y}(n)\right\Vert \geq \frac{\varepsilon y}{2n^{2}}%
\right) ,  \label{IMM}
\end{eqnarray}%
where the last equality is justified by the fact that the event $\left\{ \mu
<n<\Upsilon \right\} $ is specified by $\mathbf{V}(0),\mathbf{V}(1),...,\mathbf{V}(n-1)$ while $\mathbf{%
Y}(n)$ depends on the environment and immigration after the moment
$n-1$. Besides, the manipulations with sums are correct since the
sums run only
till $\Upsilon $ which is finite with probability 1 by Lemma \ref%
{LtailAsympt}.

Clearly,%
\begin{equation}
\mathbf{P}\left( \left\Vert \mathbf{Y}(n)\right\Vert \geq \frac{\varepsilon y%
}{2n^{2}}\right) =\mathbf{P}\left( \left\Vert \mathbf{Y}(1)\right\Vert \geq
\frac{\varepsilon y}{2n^{2}}\right) .  \label{IIMM}
\end{equation}%
Now using Lemma \ref{LEstzero} we have
\begin{equation*}
\mathbf{P}\left( \left\Vert \mathbf{Y}(1)\right\Vert \geq \frac{\varepsilon y%
}{2n^{2}}\right) \leq \frac{2^{\kappa }Kn^{2\kappa }}{\varepsilon ^{\kappa
}y^{\kappa }}.
\end{equation*}%
Thus,
\begin{eqnarray*}
\mathbf{P}\left( \sum_{n=\mu +1}^{\Upsilon -1}\left\Vert \mathbf{Y}%
(n)\right\Vert \geq \varepsilon y,\mu <\Upsilon \right) &\leq& \frac{K_{1}}{%
\varepsilon ^{\kappa }y^{\kappa }}\sum_{n=1}^{\infty }n^{2\kappa }\mathbf{P}%
\left( \mu <n<\Upsilon \right) \\
&\leq &\frac{K_{2}}{\varepsilon ^{\kappa }y^{\kappa }}\mathbf{E}\left[
\Upsilon ^{2\kappa +1}I\left\{ \mu <\Upsilon \right\} \right] \leq \frac{%
\varepsilon }{y^{\kappa }}
\end{eqnarray*}%
for $r\geq r_{1}(\varepsilon )$ since $\mathbf{E}\Upsilon ^{2\kappa
+1}<\infty $ in view of Lemma \ref{LtailAsympt} and $\mu =\mu (r)\uparrow
\infty $ in probability as $r\rightarrow \infty .$

The lemma is proved.

\begin{lemma}
\label{Lim2}If \
\begin{equation*}
\mathbf{E}\left( \left\Vert A\right\Vert +\mathbf{E}\left\Vert \mathbf{\eta }%
\right\Vert \right) ^{\kappa }<\infty
\end{equation*}%
then \ for any fixed $r$%
\begin{equation*}
\mathbf{E}\left[ \left\Vert \mathbf{V}(\mu )\right\Vert ^{\kappa }I\left\{
\mu <\Upsilon \right\} \right] <\infty .
\end{equation*}
\end{lemma}

\textbf{Proof}. It is easy to check that for any $y,w\geq 0,$ any $%
\varepsilon \in (0,1)$ and $x>0$%
\begin{equation*}
y^{x}\leq (1+\varepsilon )w^{x}+K\left\vert y-w\right\vert ^{x},
\end{equation*}%
where $K=K(x,\varepsilon ):=\left( 1-\left( 1+\varepsilon \right)
^{-1/x}\right) ^{-x}$. Hence introducing a temporary  notation $\mathbf{V}%
[0,n):=\left( \mathbf{V}(0),..,\mathbf{V}(n-1)\right) $ we have
\begin{eqnarray}
&&\mathbf{E}\left[ \left\Vert \mathbf{V}(n)\right\Vert ^{x}|\mathbf{V}(0),..,%
\mathbf{V}(n-1)\right]\qquad\qquad\qquad\qquad\qquad  \notag \\
&&  \notag \\
&&\quad=\mathbf{E}\left[ \left\Vert \mathbf{V}(n-1)A_{n-1}+\mathbf{V}(n)-\mathbf{V%
}(n-1)A_{n-1}\right\Vert ^{x}|\mathbf{V}[0,n)\right]  \notag \\
&&  \notag \\
&&\quad\leq (1+\varepsilon )\left\Vert \mathbf{V}(n-1)A_{n-1}\right\Vert ^{x}+K%
\mathbf{E}\left[ \left\Vert
\mathbf{V}(n)-\mathbf{V}(n-1)A_{n-1}\right\Vert
^{x}|\mathbf{V}[0,n)\right]  \notag \\
&&  \notag \\
&&\quad\leq (1+\varepsilon )\left\Vert
\mathbf{V}(n-1)A_{n-1}\right\Vert
^{x}+2^{x}K\mathbf{E}\left\Vert \mathbf{\eta }(n)\right\Vert ^{x}  \notag \\
&&  \notag \\
&&\,\quad+2^{x}K\sum_{i=1}^{m}\mathbf{E}\left[ \left\Vert
\sum_{k=1}^{V_{i}(n-1)}\sum_{j=1}^{m}\left[ \mathbf{\xi }%
_{ij}(n-1;k)-a_{ij}(n-1)\right] \right\Vert
^{x}\Big|\,\mathbf{V}[0,n)\right]. \label{SSta}
\end{eqnarray}%
Recalling the definition $a_{ij}(n)=\mathbf{E}_{\mathbf{f}_{n}}\xi _{ij}(n)$%
, set%
\begin{equation}
\beta _{i}(n):=\sum_{j=1}^{m}(\mathbf{\xi }_{ij}(n)-a_{ij}(n))
\label{DefetaBeta}
\end{equation}%
and let%
\begin{equation*}
M_{x}(n;i):=\mathbf{E}_{\mathbf{f}}\left\vert \beta _{i}(n)\right\vert
^{x},\ \ M_{x}:=\max_{1\leq i\leq m}\mathbf{E}\left\vert \beta
_{i}(n)\right\vert ^{x}=\max_{1\leq i\leq m}\mathbf{E}\left\vert
\sum_{j=1}^{m}(\mathbf{\xi }_{ij}-a_{ij})\right\vert ^{x}.
\end{equation*}%
Observe that $M_{\kappa }<\infty $ in view of condition (\ref{Ckappabig}).

By Lemma \ref{LGut} for $x>1$ we have
\begin{eqnarray*}
&&\mathbf{E}\left[ \left\Vert \sum_{k=1}^{V_{i}(n-1)}\sum_{j=1}^{m}\left[
\mathbf{\xi }_{ij}(n-1;k)-a_{ij}(n-1)\right] \right\Vert ^{x}\Big|\,\mathbf{V}[0,n)%
\right] \\
&&\quad\leq m^{x}\sum_{j=1}^{m}\mathbf{E}\left[ \left\Vert \sum_{k=1}^{V_{i}(n-1)}%
\left[ \mathbf{\xi }_{ij}(n-1;k)-a_{ij}(n-1)\right] \right\Vert ^{x}\Big|\,\mathbf{%
V}(n-1)\right] \\
&&\quad\quad\leq
R_{x}m^{x}\sum_{i=1}^{m}M_{x}(n-1,i)V_{i}^{x/2\vee 1}(n-1).
\end{eqnarray*}%
Thus,%
\begin{eqnarray}
\mathbf{E}\left[ \left\Vert \mathbf{V}(n)\right\Vert ^{x}|\mathbf{V}[0,n)%
\right] &\leq &(1+\varepsilon )\left\Vert \mathbf{V}(n-1)A_{n-1}\right\Vert
^{x}+2^{x}K\mathbf{E}\left\Vert \mathbf{\eta }(n)\right\Vert ^{x}  \notag \\
&&+2^{x}KR_{x}m^{x}\sum_{i=1}^{m}M_{x}(n-1,i)V_{i}^{x/2\vee 1}(n-1)  \notag
\\
&&\qquad=:\Lambda _{x}(n-1).  \label{More1}
\end{eqnarray}%
On the other hand, for $x\leq 1$%
\begin{eqnarray}
\mathbf{E}\left[ \left\Vert \mathbf{V}(n)\right\Vert ^{x}|\mathbf{V}[0,n)%
\right] &\leq &\mathbf{E}\left[ \left\Vert \mathbf{V}(n)-\mathbf{\eta }%
(n)\right\Vert ^{x}+\left\Vert \mathbf{\eta }(n)\right\Vert ^{x}|\mathbf{V}%
[0,n)\right]  \notag \\
&\leq &\left( \mathbf{E}\left[ \left\Vert \mathbf{V}(n)-\mathbf{\eta }%
(n)\right\Vert |\mathbf{V}[0,n)\right] \right) ^{x}+\mathbf{E}\left\Vert
\mathbf{\eta }\right\Vert ^{x}  \notag \\
&\leq &\left\Vert \mathbf{V}(n-1)A_{n-1}\right\Vert ^{x}+\mathbf{E}%
\left\Vert \mathbf{\eta }\right\Vert ^{x}\leq \Lambda _{x}(n-1).
\label{Less1}
\end{eqnarray}%
Clearly,%
\begin{eqnarray*}
\Lambda _{\kappa }(n-1)I\left\{ n-1<\mu <\Upsilon \right\} &\leq
&Q_{n-1}(r):=(1+\varepsilon )r^{\kappa }\left\Vert A_{n-1}\right\Vert
^{\kappa }+2^{\kappa }K\mathbf{E}\left\Vert \mathbf{\eta }\right\Vert
^{\kappa } \\
&&+2^{\kappa }KR_{\kappa }m^{\kappa }r^{\kappa /2\vee
1}\sum_{i=1}^{m}M_{\kappa }(n-1,i).
\end{eqnarray*}%
Observe that by our conditions the expectation $\mathbf{E}\left[ Q_{n-1}(r)%
\right] $ is independent of $n$ and, moreover,%
\begin{eqnarray*}
\mathbf{E}\left[ Q_{n-1}(r)\right] &=&(1+\varepsilon )r^{\kappa }\mathbf{E}%
\left\Vert A\right\Vert ^{\kappa }+2^{\kappa }K\mathbf{E}\left\Vert \mathbf{%
\eta }\right\Vert ^{\kappa } \\
&&+2^{\kappa }KR_{\kappa }m^{\kappa }r^{\kappa /2\vee 1}\mathbf{E}\left[
\sum_{i=1}^{m}M_{\kappa }(n-1,i)\right] \\
&\leq &(1+\varepsilon )r^{\kappa }\mathbf{E}\left\Vert A\right\Vert ^{\kappa
}+2^{\kappa }K\mathbf{E}\left\Vert \mathbf{\eta }\right\Vert ^{\kappa
}+2^{\kappa }KR_{\kappa }m^{\kappa +1}r^{\kappa /2\vee 1}M_{\kappa
}:=K_{1}<\infty .
\end{eqnarray*}
Besides,
\begin{eqnarray*}
\left\Vert \mathbf{V}(\mu )\right\Vert ^{\kappa }I\left\{ \mu <\Upsilon
\right\} &=&\Lambda _{\kappa }(\mu -1)\frac{\left\Vert \mathbf{V}(\mu
)\right\Vert ^{\kappa }}{\Lambda _{\kappa }(\mu -1)}I\left\{ \mu <\Upsilon
\right\} \leq Q_{\mu -1}(r)\frac{\left\Vert \mathbf{V}(\mu )\right\Vert
^{\kappa }}{\Lambda _{\kappa }(\mu -1)} \\
&\leq &\sum_{\mu \leq n<\Upsilon }Q_{n-1}(r)\frac{\left\Vert \mathbf{V}%
(n)\right\Vert ^{\kappa }}{\Lambda _{\kappa }(n-1)}.
\end{eqnarray*}%
Now by (\ref{More1}) and (\ref{Less1}) we see that
\begin{eqnarray*}
&&\mathbf{E}\left[ \left\Vert \mathbf{V}(\mu )\right\Vert ^{\kappa }I\left\{
\mu <\Upsilon \right\} \right] \leq \mathbf{E}\left[ \sum_{\mu \leq
n<\Upsilon }Q_{n-1}(r)\frac{\left\Vert \mathbf{V}(n)\right\Vert ^{\kappa }}{%
\Lambda _{\kappa }(n-1)}\right] \\
&&\quad\leq\mathbf{E}\left[ \sum_{n=1}^{\infty
}Q_{n-1}(r)\frac{\left\Vert \mathbf{V}(n)\right\Vert ^{\kappa
}}{\Lambda _{\kappa }(n-1)}I\left\{
\Upsilon \geq n\right\} \right] \\
&&\quad=\sum_{n=1}^{\infty }\mathbf{E}\left[
\frac{Q_{n-1}(r)}{\Lambda _{\kappa
}(n-1)}I\left\{ \Upsilon \geq n\right\} \mathbf{E}\left[ \left\Vert \mathbf{V%
}(n)\right\Vert ^{\kappa }|\mathbf{V}[0,n)\right] \right] \\
&&\quad\quad\leq \sum_{n=1}^{\infty }\mathbf{E}\left[
Q_{n-1}(r)I\left\{ \Upsilon \geq n\right\} \right]
=\sum_{n=1}^{\infty }\mathbf{E}\left[ Q_{n-1}(r)\right]
\mathbf{P}\left( \Upsilon \geq n\right) \\
&&\qquad\qquad\qquad\qquad\qquad\qquad\qquad\quad\leq
K_{1}\sum_{n=1}^{\infty }\mathbf{P}\left( \Upsilon \geq n\right)
<\infty .
\end{eqnarray*}%
The lemma is proved.

Let \
\begin{equation*}
\mathbf{Z}(i,j;l,n):=\left( Z_{1}(i,j;l,n),...,Z_{1}(i,j;l,n)\right)
\end{equation*}%
be the $m-$dimensional vector of the total progeny alive at moment $n>l$ of
the $j-$th particle of type $i$ existing in the process at time $l$. Let,
further, $\mathbf{Z}\left( l,l\right) :=\mathbf{V}(l)$ and
\begin{equation*}
\mathbf{Z}\left( l,n\right) =\left(
Z_{1}(l,n),...,Z_{m}(l,n)\right)
:=\sum_{i=1}^{m}\sum_{j=1}^{V_{i}(l)}\mathbf{Z}(i,j;l,n)
\end{equation*}%
be the total progeny at time $n$ of all particles existing in the process at
moment $l<n$. Recalling definition (\ref{finprod}) we set
\begin{equation*}
\Gamma \left( \mu ,\mu \right) :=\left\langle \mathbf{V}(\mu ),\mathbf{C}%
_{\mu }\right\rangle =\left\langle \mathbf{Z}(\mu ,\mu ),\mathbf{C}_{\mu
}\right\rangle
\end{equation*}%
and for $n>\mu $ let
\begin{equation}
\Gamma \left( \mu ,n\right) :=\sum_{i=1}^{m}\sum_{j=1}^{V_{i}(\mu
)}\left\langle \mathbf{Z}(i,j;\mu ,n),\mathbf{C}_{n}\right\rangle
=\left\langle \mathbf{Z}(\mu ,n),\mathbf{C}_{n}\right\rangle
\label{DefGamMUMU}
\end{equation}%
be the inner product of $\mathbf{C}_{n}$ and the vector of the
total progeny at time $n$ of the $\mathbf{V}(\mu )$ particles
presented in the process at time $\mu,$   and
\begin{equation}
\Gamma (\mu ):=\sum_{n=\mu }^{\Upsilon -1}\Gamma \left( \mu ,n\right)
=\sum_{n=\mu }^{\Upsilon -1}\left\langle \mathbf{Z}(\mu ,n),\mathbf{C}%
_{n}\right\rangle .  \label{DefGammaMu}
\end{equation}

Denote by $\mathcal{C}_{n}$ ($n=1,2,\ldots $) be the $\sigma $-algebra
generated by the tuple%
\begin{equation*}
\mathbf{H}_{0}(\mathbf{s},\lambda ),\mathbf{H}_{1}(\mathbf{s},\lambda
),\ldots ,\mathbf{H}_{n-1}(\mathbf{s},\lambda ),\mathbf{V}(0),\ldots ,%
\mathbf{V}(n).
\end{equation*}%
The next lemma shows that for large $r$ the random variable $\Gamma \left(
\mu \right) =\Gamma \left( \mu (r)\right) $ is, in a sense, close to the
conditional expectation $\mathbf{E}\left[ \Gamma \left( \mu \right) |%
\mathcal{C}_{\mu }\right] =\mathbf{V}(\mu )\Xi _{\mu }$ (recall definition (%
\ref{DefSumThet})).

\begin{lemma}
\label{Lke2 2}Under the conditions of Theorem \ref{TTmain} for any $%
\varepsilon >0$ there exists $r_{1}=r_{1}(\varepsilon )$ such that for all $%
y>0$%
\begin{equation}
\mathbf{P}\left( \left\vert \Gamma \left( \mu \right) -\left\langle \mathbf{V%
}(\mu ),\Xi _{\mu }\right\rangle \right\vert >\varepsilon y;\mu <\Upsilon
\right) \leq \frac{\varepsilon }{y^{\kappa }}\mathbf{E}\left[ \left\Vert
\mathbf{V}(\mu )\right\Vert ^{\kappa }I\left\{ \mu <\Upsilon \right\} \right]
.  \label{Dif11}
\end{equation}
\end{lemma}

\textbf{Proof}. Evidently, for $n\geq \mu +1$%
\begin{eqnarray*}
\left\langle \mathbf{Z}(\mu ,n)-\mathbf{V}(\mu )\Pi _{\mu ,n},\mathbf{C}%
_{n}\right\rangle &=&\sum_{l=\mu +1}^{n}\left\langle \mathbf{Z}(\mu ,l)\Pi
_{l,n}-\mathbf{Z}(\mu ,l-1)\Pi _{l-1,n},\mathbf{C}_{n}\right\rangle \\
&=&\sum_{l=\mu +1}^{n}\left\langle \mathbf{Z}(\mu ,l)-\mathbf{Z}(\mu
,l-1)A_{l-1},\Pi _{l,n}\mathbf{C}_{n}\right\rangle
\end{eqnarray*}%
which implies%
\begin{eqnarray}
\left\vert \Gamma \left( \mu \right) -\left\langle \mathbf{V}(\mu ),\Xi
_{\mu }\right\rangle \right\vert &=&\left\vert \sum_{n=\mu +1}^{\infty
}\left\langle \mathbf{Z}(\mu ,n)-\mathbf{V}(\mu )\Pi _{\mu ,n},\mathbf{C}%
_{n}\right\rangle \right\vert  \notag \\
&\leq &\sum_{n=\mu +1}^{\infty }\left\vert \left\langle \mathbf{Z}(\mu ,n)-%
\mathbf{V}(\mu )\Pi _{\mu ,n},\mathbf{C}_{n}\right\rangle \right\vert  \notag
\\
&=&\sum_{n=\mu +1}^{\infty }\left\vert \sum_{l=\mu +1}^{n}\left\langle
\mathbf{Z}(\mu ,l)-\mathbf{Z}(\mu ,l-1)A_{l-1},\Pi _{l,n}\mathbf{C}%
_{n}\right\rangle \right\vert  \notag \\
&\leq &\sum_{l=\mu +1}^{\infty }\sum_{n=l}^{\infty }\left\vert \left\langle
\mathbf{Z}(\mu ,l)-\mathbf{Z}(\mu ,l-1)A_{l-1},\Pi _{l,n}\mathbf{C}%
_{n}\right\rangle \right\vert  \notag \\
&\leq &\sum_{l=\mu +1}^{\infty }\left\Vert \mathbf{Z}(\mu ,l)-\mathbf{Z}(\mu
,l-1)A_{l-1}\right\Vert \left\Vert \sum_{n=l}^{\infty }\Pi _{l,n}\mathbf{C}%
_{n}\right\Vert  \notag \\
&=&\sum_{l=\mu +1}^{\infty }\left\Vert \mathbf{Z}(\mu ,l)-\mathbf{Z}(\mu
,l-1)A_{l-1}\right\Vert \left\Vert \Xi _{l}\right\Vert .  \label{seco2}
\end{eqnarray}%
Here the interchange of summation order is justified since the
process having $\mathbf{V}(\mu )$ particles at moment $\mu
<\Upsilon <\infty $ dies out with probability 1.

Thus,%
\begin{eqnarray*}
&&\mathbf{P}\left( \left\vert \Gamma \left( \mu \right) -\left\langle
\mathbf{V}(\mu ),\Xi _{\mu }\right\rangle \right\vert >\varepsilon y;\mu
<\Upsilon \right) \\
&&\qquad\leq \mathbf{P}\left( \sum_{l=\mu +1}^{\infty }\left\Vert
\mathbf{Z}(\mu ,l)-\mathbf{Z}(\mu ,l-1)A_{l-1}\right\Vert
\left\Vert \Xi _{l}\right\Vert
>\varepsilon y;\mu <\Upsilon \right) .
\end{eqnarray*}%
Hence, using (\ref{Series}) we conclude%
\begin{eqnarray*}
&&\mathbf{P}\left( \left\vert \Gamma \left( \mu \right) -\left\langle
\mathbf{V}(\mu ),\Xi _{\mu }\right\rangle \right\vert >\varepsilon y;\mu
<\Upsilon |\mathcal{C}_{\mu }\right) \\
&\leq &\mathbf{P}\left( \sum_{n=1}^{\infty }\left\Vert \mathbf{Z}(\mu ,n+\mu
)-\mathbf{Z}(\mu ,n+\mu -1)A_{n+\mu -1}\right\Vert \left\Vert \Xi _{n+\mu
}\right\Vert >\frac{\varepsilon y}{2}\sum_{n=1}^{\infty }\frac{1}{n^{2}};\mu
<\Upsilon \Big|\,\mathcal{C}_{\mu }\right) \\
&\leq &\sum_{n=1}^{\infty }\mathbf{P}\left( \left\Vert \mathbf{Z}(\mu ,n+\mu
)-\mathbf{Z}(\mu ,n+\mu -1)A_{n+\mu -1}\right\Vert \left\Vert \Xi _{n+\mu
}\right\Vert >\frac{\varepsilon y}{2}\frac{1}{n^{2}};\mu <\Upsilon\Big|\,\mathcal{%
C}_{\mu }\right) .
\end{eqnarray*}%
Since $\mathbf{Z}(\mu ,n+\mu )-\mathbf{Z}(\mu ,n+\mu -1)A_{n+\mu -1}$ and $%
\Xi _{n+\mu }$ are independent random objects on the event $\mu =k<\Upsilon
<\infty $ and $\Xi _{n+\mu }\overset{d}{=}\Xi $ (see (\ref{DefSumThet})), we
get%
\begin{eqnarray*}
&&\mathbf{P}\left( \left\Vert \mathbf{Z}(\mu ,n+\mu )-\mathbf{Z}(\mu ,n+\mu
-1)A_{n+\mu -1}\right\Vert \left\Vert \Xi _{n+\mu }\right\Vert >\frac{%
\varepsilon y}{2}\frac{1}{n^{2}};\mu <\Upsilon \Big|\,\mathcal{C}_{\mu }\right) \\
&&\quad =\int_{0}^{\infty }\mathbf{P}\left( \left\Vert
\mathbf{Z}(\mu ,n+\mu )-\mathbf{Z}(\mu ,n+\mu -1)A_{n+\mu
-1}\right\Vert \in dt\,;\,\mu <\Upsilon |\mathcal{C}_{\mu }\right)
\mathbf{P}\left( \left\Vert \Xi \right\Vert \geq \frac{\varepsilon
y}{2tn^{2}}\right) .
\end{eqnarray*}%
According to Condition $T$ there exists a constant $K_{1}\in \left( 0,\infty
\right) $ such that for all $n>0$
\begin{eqnarray}
&&\mathbf{P}\left( \left\Vert \mathbf{Z}(\mu ,n+\mu )-\mathbf{Z}(\mu ,n+\mu
-1)A_{n+\mu -1}\right\Vert \left\Vert \Xi _{n+\mu }\right\Vert >\frac{%
\varepsilon y}{2}\frac{1}{n^{2}};\mu
<\Upsilon\Big|\,\mathcal{C}_{\mu }\right)
\notag \\
&&\quad \leq \int_{0}^{\infty }\mathbf{P}\left( \left\Vert \mathbf{Z}(\mu
,n+\mu )-\mathbf{Z}(\mu ,n+\mu -1)A_{n+\mu -1}\right\Vert \in dt\,;\,\mu
<\Upsilon |\mathcal{C}_{\mu }\right) \frac{K_{1}t^{\kappa }}{\varepsilon
^{\kappa }y^{\kappa }}n^{2\kappa }  \notag \\
&&\quad \leq \frac{K_{1}}{\varepsilon ^{\kappa }y^{\kappa }}n^{2\kappa }%
\mathbf{E}\left[ \left\Vert \mathbf{Z}(\mu ,n+\mu )-\mathbf{Z}(\mu ,n+\mu
-1)A_{n+\mu -1}\right\Vert ^{\kappa }I\left\{ \mu <\Upsilon \right\}
\left\vert \,\mathcal{C}_{\mu }\right. \right] .  \label{BB}
\end{eqnarray}

Now we consider the cases $\kappa \leq 1$ and $\kappa >1$ separately.

For the first case we use for $l>\mu $ \ the estimate%
\begin{eqnarray}
&&\mathbf{E}\left[ \left\Vert \mathbf{Z}(\mu ,l)-\mathbf{Z}(\mu
,l-1)A_{l-1}\right\Vert ^{\kappa }I\left\{ \mu <\Upsilon \right\} \left\vert
\,\mathcal{C}_{\mu }\right. \right]  \notag \\
&&\qquad\leq \left( \mathbf{E}\left[ \left\Vert \mathbf{Z}(\mu
,l)-\mathbf{Z}(\mu
,l-1)A_{l-1}\right\Vert ^{2}I\left\{ \mu <\Upsilon \right\} \left\vert \,%
\mathcal{C}_{\mu }\right. \right] \right) ^{\kappa /2}.  \label{Nstarr}
\end{eqnarray}%
Further, we have (recall (\ref{DefetaBeta}))
\begin{eqnarray}
&&\mathbf{E}\left[ \left\Vert \mathbf{Z}(\mu ,l)-\mathbf{Z}(\mu
,l-1)A_{l-1}\right\Vert ^{2}I\left\{ \mu <\Upsilon \right\} \left\vert \,%
\mathcal{C}_{\mu }\right. \right]  \notag \\
&&\quad =\mathbf{E}\left[ \left( \sum_{i=1}^{m}\sum_{k=1}^{Z_{i}(\mu
,l-1)}\sum_{j=1}^{m}\left[ \mathbf{\xi }_{ij}(l-1;k)-a_{ij}(l-1)\right]
\right) ^{2}I\left\{ \mu <\Upsilon \right\} \left\vert \,\mathcal{C}_{\mu
}\right. \right]  \notag \\
&&\quad \,=\sum_{i=1}^{m}\mathbf{E}_{\mathbf{f}}\beta _{i}^{2}(l-1)\mathbf{E}%
\left[ Z_{i}(\mu ,l-1)I\left\{ \mu <\Upsilon \right\} \left\vert \,\mathcal{C%
}_{\mu }\right. \right]  \notag \\
&&\qquad=\sum_{i=1}^{m}\mathbf{E}_{\mathbf{f}}\beta _{i}^{2}(l-1)\left( \mathbf{V}%
(\mu )\Pi _{\mu ,l-1}\right) _{i}\leq \left\Vert \mathbf{V}(\mu
)\Pi _{\mu ,l-1}\right\Vert
\sum_{i=1}^{m}\mathbf{E}_{\mathbf{f}}\beta _{i}^{2}(l-1).
\notag \\
&&{}  \label{BB11}
\end{eqnarray}%
Thus, for $\kappa \leq 1$ and $l=n+\mu >\mu $
\begin{eqnarray}
&&\left( \mathbf{E}\left[ \left\Vert \mathbf{Z}(\mu ,n+\mu )-\mathbf{Z}(\mu
,n+\mu -1)A_{n+\mu -1}\right\Vert ^{2}I\left\{ \mu <\Upsilon \right\}
\left\vert \,\mathcal{C}_{\mu }\right. \right] \right) ^{\kappa /2}  \notag
\\
& &\qquad\qquad\qquad\leq\left\Vert \mathbf{V}(\mu )\Pi _{\mu
,n+\mu -1}\right\Vert ^{\kappa /2}\left(
\sum_{i=1}^{m}\mathbf{E}_{\mathbf{f}}\beta _{i}^{2}(n+\mu
-1)\right) ^{\kappa /2}.  \label{NNstar}
\end{eqnarray}%
This, in view of the inequality
\begin{equation*}
s(\kappa /2)=\lim_{n\rightarrow \infty }\left( \mathbf{E}\left\Vert \Pi
_{0,n}\right\Vert ^{\kappa /2}\right) ^{1/n}<1,
\end{equation*}%
the first part of condition (\ref{Ckappasmall}), and relations (\ref{BB11})-(%
\ref{NNstar}) leads to the estimate%
\begin{eqnarray*}
&&\mathbf{P}\left( \left\vert \Gamma \left( \mu \right) -\left\langle
\mathbf{V}(\mu ),\Xi _{\mu }\right\rangle \right\vert >\varepsilon y;\mu
<\Upsilon \right) \\
& &\leq\frac{K}{\varepsilon ^{\kappa }y^{\kappa }}\mathbf{E}\left[
\sum_{n=1}^{\infty }n^{2\kappa }\left\Vert \mathbf{V}(\mu
)\right\Vert ^{\kappa /2}\left\Vert \Pi _{\mu ,n+\mu
-1}\right\Vert ^{\kappa /2}I\left\{
\mu <\Upsilon <\infty \right\} \right] \\
&&=\frac{K}{\varepsilon ^{\kappa }y^{\kappa }}\mathbf{E}\left[
\left\Vert \mathbf{V}(\mu )\right\Vert ^{\kappa /2}I\left\{ \mu
<\Upsilon <\infty \right\} \sum_{n=1}^{\infty }n^{2\kappa
}\mathbf{E}\left\Vert \Pi
_{0,n-1}\right\Vert ^{\kappa /2}\right] \\
&&\leq \frac{K_{1}}{\varepsilon ^{\kappa }y^{\kappa }r^{\kappa /2}}\mathbf{E}%
\left[ \left\Vert \mathbf{V}(\mu )\right\Vert ^{\kappa }I\left\{
\mu
<\Upsilon <\infty \right\} \right] \leq \frac{\varepsilon }{y^{\kappa }}%
\mathbf{E}\left[ \left\Vert \mathbf{V}(\mu )\right\Vert ^{\kappa }I\left\{
\mu <\Upsilon <\infty \right\} \right]
\end{eqnarray*}%
for all $r\geq r_{0}(\varepsilon )$, proving the lemma for $\kappa \leq 1$.

For the case $\kappa >1$ we use Lemma \ref{LGut} to conclude that for any $%
l>\mu $%
\begin{eqnarray*}
&&\mathbf{E}\left[ \left\Vert \mathbf{Z}(\mu ,l)-\mathbf{Z}(\mu
,l-1)A_{l-1}\right\Vert ^{\kappa }I\left\{ \mu <\Upsilon <\infty \right\}
\left\vert \mathcal{C}_{\mu }\right. \right] \\
&&\quad \leq R_{\kappa }m^{\kappa }\sum_{i=1}^{m}M_{\kappa }(l;i)\mathbf{E}%
\left[ \left\Vert \mathbf{Z}\left( \mu ,l-1\right) \right\Vert ^{\kappa
/2\vee 1}I\left\{ \mu <\Upsilon <\infty \right\} \left\vert \mathcal{C}_{\mu
}\right. \right] \\
&&\quad \leq R_{\kappa }m^{\kappa }\mathbf{E}\left[ \left\Vert \mathbf{Z}%
\left( \mu ,l-1\right) \right\Vert ^{\kappa /2\vee 1}I\left\{ \mu <\Upsilon
<\infty \right\} \left\vert \mathcal{C}_{\mu }\right. \right]
\sum_{i=1}^{m}M_{\kappa }(l-1;i).
\end{eqnarray*}%
By Lemma \ref{LGeom} there exist constants $K=K\left( \kappa /2\vee 1\right)
<\infty $ and $\rho =\rho (\kappa /2\vee 1)\in \left( 0,1\right) $ such that
\begin{equation*}
\mathbf{E}\left[ \left\Vert \mathbf{Z}\left( \mu ,l-1\right) \right\Vert
^{\kappa /2\vee 1}I\left\{ \mu <\Upsilon <\infty \right\} \left\vert
\mathcal{C}_{\mu }\right. \right] \leq K\rho ^{l-\mu -1}\left\Vert \mathbf{V}%
(\mu )\right\Vert ^{\kappa /2\vee 1}.
\end{equation*}%
This inequality combined with condition (\ref{Ckappabig}) yields for $%
l=n+\mu $ the estimates
\begin{eqnarray*}
&&\mathbf{P}\left( \left\vert \Gamma \left( \mu \right) -\left\langle
\mathbf{V}(\mu ),\Xi _{\mu }\right\rangle \right\vert >\varepsilon y;\mu
<\Upsilon \right) \\
&\leq &\frac{K_{1}R_{\kappa }m^{\kappa }}{\varepsilon ^{\kappa }y^{\kappa }}%
\mathbf{E}\left[ \sum_{n=1}^{\infty }n^{2\kappa }\mathbf{E}\left[
\left\Vert \mathbf{Z}\left( \mu ,n+\mu -1\right) \right\Vert
^{\kappa /2\vee 1}\sum_{i=1}^{m}M_{\kappa }(n+\mu -1;i)\big|
\mathcal{C}_{\mu }
\right] I\left\{ \mu <\Upsilon <\infty \right\} \right] \\
&\leq &\frac{K_{1}R_{\kappa }m^{\kappa +1}M_{\kappa }}{\varepsilon ^{\kappa
}y^{\kappa }}\mathbf{E}\left[ \sum_{n=1}^{\infty }n^{2\kappa }\mathbf{E}%
\left[ \left\Vert \mathbf{Z}\left( \mu ,n+\mu -1\right) \right\Vert ^{\kappa
/2\vee 1}\left\vert \,\mathcal{C}_{\mu }\right. \right] I\left\{ \mu
<\Upsilon <\infty \right\} \right] \\
&\leq &\frac{K_{1}R_{\kappa }m^{\kappa +1}M_{\kappa }K}{\varepsilon ^{\kappa
}y^{\kappa }}\mathbf{E}\left[ \sum_{n=1}^{\infty }n^{2\kappa }\left\Vert
\mathbf{V}(\mu )\right\Vert ^{\kappa /2\vee 1}\rho ^{n-1}I\left\{ \mu
<\Upsilon <\infty \right\} \right] \\
&=&\frac{K_{2}}{\varepsilon ^{\kappa }y^{\kappa }}\mathbf{E}\left[
\left\Vert \mathbf{V}(\mu )\right\Vert ^{\kappa /2\vee 1}I\left\{ \mu
<\Upsilon \right\} \sum_{n=1}^{\infty }n^{2\kappa }\rho ^{n-1}\right] \\
&\leq &\frac{K_{3}}{\varepsilon ^{\kappa }y^{\kappa }r^{\kappa -\kappa
/2\vee 1}}\mathbf{E}\left[ \left\Vert \mathbf{V}(\mu )\right\Vert ^{\kappa
}I\left\{ \mu <\Upsilon \right\} \right] \leq \frac{\varepsilon }{y^{\kappa }%
}\mathbf{E}\left[ \left\Vert \mathbf{V}(\mu )\right\Vert ^{\kappa }I\left\{
\mu <\Upsilon \right\} \right]
\end{eqnarray*}%
(the last is valid by selecting $r$ sufficiently large) which justifies the
statement of the lemma for $\kappa >1.$

The lemma is proved.

\section{The accumulated amount of the final product}

In this section we deduce some estimates related with the total size of the
final product accumulated in a subcritical MBPIFPRE within a life-period.

Let
\begin{equation*}
\Delta _{V}(0):=\Theta (0)+\psi (1)I\left\{ \left\Vert \mathbf{\eta }%
(1)\right\Vert >0\right\}
\end{equation*}%
and let

\begin{equation*}
\Delta _{V}(n):=\sum_{i=1}^{m}\sum_{k=1}^{V_{i}(n)}\varphi _{i}(n;k),n\geq 1,
\end{equation*}%
be the total amount of the final product produced by the particles of the $n$%
-th generation of a MBPIFPRE,%
\begin{equation*}
\Delta _{I}(n)=\Delta _{V}(n)+\psi \left( n+1\right)
\end{equation*}%
be the total amount of the final product produced by the individuals of the $%
n$-th generation of a MBPIFPRE plus the final product contributed by the
particles immigrating at moment $n+1$, and
\begin{equation*}
Y_{I}(n):=\psi \left( n+1\right) +\sum_{i=1}^{m}\sum_{j=1}^{\eta
_{i}(n+1)}\left( \varphi _{i}(n+1,j)+\sum_{k=n+2}^{\Upsilon
-1}\sum_{t=1}^{m}\sum_{w=1}^{Z_{t}^{I}(i,j;n+1,k)}\varphi _{t}(k,w)\right)
\end{equation*}%
be the total size of the final product contributed to the system
by all particles immigrating in the process at moment $n+1$ and by
their progeny.

Of interest will be the quantity
\begin{equation*}
\Theta _{Z}(\mu ,n):=\sum_{t=1}^{m}\sum_{k=1}^{Z_{t}(\mu ,n)}\varphi
_{t}(n;k),
\end{equation*}%
i.e., the total amount of the final product produced by the particles of the
$n$-th generation, $n>\mu ,$ of the MBPIFPRE which belong to the progeny of
the $\mathbf{V}(\mu )$ particles of generation $\mu $ and the random
variable
\begin{equation*}
\Theta _{Z}\left( \mu \right) :=\sum_{n=\mu }^{\Upsilon -1}\Theta _{Z}\left(
\mu ,n\right)
\end{equation*}%
which is equal to the total amount of the final product produced by the $%
\mathbf{V}(\mu )$ particles of generation $\mu $ and their progeny up to the
moment of extinction of the MBPRE generated by the $\mathbf{V}(\mu )$
particles of the $\mu -$th generation.

Clearly, the total amount
\begin{equation*}
\Theta :=\sum_{n=0}^{\Upsilon -1}\Delta _{I}(n)
\end{equation*}%
of the final product accumulated in the system during the life
period which
starts at moment $n=1$ admits on the set $\mu <\Upsilon $ the representation%
\begin{equation*}
\Theta =\sum_{n=0}^{\mu -1}\Delta _{I}(n)+\Theta _{Z}\left( \mu \right)
+\sum_{n=\mu }^{\Upsilon -1}Y_{I}(n).
\end{equation*}

Our aim is to investigate the asymptotic behavior of the tail
distribution of~$\Theta $.

As the first step in solving this problem we compare $\Theta
_{Z}\left( \mu \right) $ with $\Gamma (\mu )$. This will be done
by the arguments similar to those used to demonstrate Lemma
\ref{Lke2 2}.

\begin{lemma}
\label{Lke2222}Under the conditions of Theorem \ref{TTmain} for any $%
\varepsilon >0$ there exists $r=r(\varepsilon )$ such that for all $y>0$%
\begin{equation}
\mathbf{P}\left( \left\vert \Theta _{Z}\left( \mu \right) -\Gamma (\mu
)\right\vert >\varepsilon y;\mu <\Upsilon \right) \leq \frac{\varepsilon }{%
y^{\kappa }}\mathbf{E}\left[ \left\Vert \mathbf{V}(\mu )\right\Vert ^{\kappa
}I\left\{ \mu <\Upsilon \right\} \right] .  \label{Dif22}
\end{equation}
\end{lemma}

\textbf{Proof}. For $l>\mu $ we have%
\begin{equation*}
\Theta _{Z}(\mu ,l)\mathbf{-}\left\langle \mathbf{Z}(\mu ,l),\mathbf{C}%
_{l}\right\rangle \mathbf{=}\sum_{t=1}^{m}\sum_{k=1}^{\mathbf{Z}_{t}(\mu
,l)}\left( \varphi _{t}(l;k)-\mathbf{E}_{\mathbf{F}}\varphi _{t}(l)\right) .
\end{equation*}%
Now we consider separately the cases $\kappa \leq 1$ and $\kappa >1$.

For $\kappa \geq 1$ we use Lemmas \ref{LGut} and \ref{LGeom} with $\rho
=\rho (\kappa /2\vee 1)$ and $K=K\left( \kappa /2\vee 1\right) $ to get
\begin{eqnarray*}
&&\mathbf{P}\left( \left\vert \Theta _{Z}(\mu ,l)\mathbf{-}\left\langle
\mathbf{Z}(\mu ,l),\mathbf{C}_{l}\right\rangle \right\vert >\frac{%
\varepsilon y}{2n^{2}}\left\vert \mathcal{C}_{\mu }\right. \right) \\
&&\qquad\leq \frac{4^{\kappa }n^{2\kappa }}{\left( \varepsilon y\right) ^{\kappa }}%
\mathbf{E}\left[ \left\vert
\sum_{t=1}^{m}\sum_{k=1}^{\mathbf{Z}_{t}(\mu ,l)}\left( \varphi
_{t}(l;k)-\mathbf{E}_{\mathbf{f}}\varphi _{t}(l)\right)
\right\vert ^{\kappa }\left\vert \mathcal{C}_{\mu }\right. \right] \\
&&\qquad\leq (4m)^{\kappa }R_{\kappa }\frac{n^{2\kappa }}{\left(
\varepsilon y\right) ^{\kappa
}}\mathbf{E%
}\left[ \left\Vert \mathbf{Z}(\mu ,l)\right\Vert ^{\kappa /2\vee
1}\left\vert \,\mathcal{C}_{\mu }\right.
\right]\sum_{i=1}^{m}\mathbf{E}_{\mathbf{f}}\left\vert \varphi
_{t}(l)-\mathbf{E}_{\mathbf{f}}\varphi _{t}(l)\right\vert ^{\kappa } \\
&&\qquad\leq (4m)^{\kappa }R_{\kappa }K\frac{n^{2\kappa }}{\left(
\varepsilon y\right) ^{\kappa }}\left\Vert \mathbf{V}(\mu
)\right\Vert ^{\kappa /2\vee 1}\rho ^{l-\mu
}\sum_{i=1}^{m}\mathbf{E}_{\mathbf{f}}\left\vert \varphi
_{t}(l)-\mathbf{E}_{\mathbf{f}}\varphi _{t}(l)\right\vert ^{\kappa
}.
\end{eqnarray*}%
Hence, in view of condition (\ref{Ckappabig}) and with $l=\mu +n$ we have%
\begin{eqnarray*}
&&\mathbf{P}\left( \left\vert \Theta _{Z}\left( \mu \right) -\Gamma (\mu
)\right\vert >\varepsilon y;\mu <\Upsilon \right) \\
&\leq &\mathbf{E}\left[ \sum_{n=1}^{\infty }\mathbf{P}_{\mathbf{F}}\left(
\left\vert \Theta _{Z}(\mu ,\mu +n)\mathbf{-}\left\langle \mathbf{Z}(\mu
,\mu +n),\mathbf{C}_{\mu +n}\right\rangle \right\vert >\frac{\varepsilon y}{%
2n^{2}}\Big|\,\mathcal{C}_{\mu }\right) I\left( \mu <\Upsilon
\right) \right] \\
&\leq &\frac{(4m)^{\kappa }R_{\kappa }KK_{1}}{\left( \varepsilon y\right)
^{\kappa }}\mathbf{E}\left[ \left\Vert \mathbf{V}(\mu )\right\Vert ^{\kappa
/2\vee 1}\sum_{n=1}^{\infty }n^{2\kappa }\rho ^{n}I\left( \mu <\Upsilon
\right) \right] \\
&=&\frac{K_{2}}{\left( \varepsilon y\right) ^{\kappa }}\mathbf{E}\left[
\left\Vert \mathbf{V}(\mu )\right\Vert ^{\kappa /2\vee 1}I\left( \mu
<\Upsilon \right) \right] \sum_{n=1}^{\infty }n^{2\kappa }\rho ^{n} \\
&\leq &\frac{K_{3}}{\left( \varepsilon y\right) ^{\kappa }r^{\kappa -(\kappa
/2\vee 1)}}\mathbf{E}\left[ \left\Vert \mathbf{V}(\mu )\right\Vert ^{\kappa
}I\left( \mu <\Upsilon \right) \right] \leq \frac{\varepsilon }{y^{\kappa }}%
\mathbf{E}\left[ \left\Vert \mathbf{V}(\mu )\right\Vert ^{\kappa }I\left(
\mu <\Upsilon \right) \right]
\end{eqnarray*}%
for all $r\geq r(\varepsilon ).$

To analyze the case $\kappa \leq 1$ we apply for $l>\mu $ the inequality%
\begin{equation*}
\mathbf{E}\left[ \left\vert \Theta _{Z}(\mu ,l)\mathbf{-}\left\langle
\mathbf{Z}(\mu ,l),\mathbf{C}_{l}\right\rangle \right\vert ^{\kappa
}\left\vert \mathcal{C}_{\mu }\right. \right] \leq \left( \mathbf{E}\left[
\left\vert \Theta _{Z}(\mu ,l)\mathbf{-}\left\langle \mathbf{Z}(\mu ,l),%
\mathbf{C}_{l}\right\rangle \right\vert ^{2}\big|\,
\mathcal{C}_{\mu } \right] \right) ^{\kappa /2}.
\end{equation*}%
Further, we have
\begin{eqnarray*}
&&\mathbf{E}\left[ \left\vert \Theta _{Z}(\mu
,l)\mathbf{-}\left\langle \mathbf{Z}(\mu
,l),\mathbf{C}_{l}\right\rangle \right\vert ^{2}\left\vert
\mathcal{C}_{\mu }\right. \right] \leq m^{2}\sum_{i=1}^{m}\mathbf{E}_{%
\mathbf{f}}\beta _{i}^{2}(l)\mathbf{E}\left[ Z_{i}(\mu
,l)\,\left\vert
\mathcal{C}_{\mu }\right. \right] \\
&&\qquad\quad=m^{2}\sum_{i=1}^{m}\mathbf{E}_{\mathbf{f}}\beta _{i}^{2}(l)\left( \mathbf{%
V}(\mu )\Pi _{\mu ,l}\right) _{i}\leq m^{2}\left\Vert
\mathbf{V}(\mu )\Pi _{\mu ,l}\right\Vert
\sum_{i=1}^{m}\mathbf{E}_{\mathbf{f}}\beta _{i}^{2}(l).
\end{eqnarray*}%
Thus, for $\kappa \leq 1$
\begin{equation*}
\left( \mathbf{E}\left[ \left\vert \Theta _{Z}(\mu ,l)\mathbf{-}\left\langle
\mathbf{Z}(\mu ,l),\mathbf{C}_{l}\right\rangle \right\vert ^{2}\left\vert
\mathcal{C}_{\mu }\right. \right] \right) ^{\kappa /2}\leq m^{\kappa
}\left\Vert \mathbf{V}(\mu )\Pi _{\mu ,l}\right\Vert ^{\kappa /2}\left(
\sum_{i=1}^{m}\mathbf{E}_{\mathbf{f}}\beta _{i}^{2}(l)\right) ^{\kappa /2}.
\end{equation*}

This combined with the assumption (\ref{Ckappasmall}) shows
that for $l=\mu +n$%
\begin{eqnarray*}
&&\mathbf{P}\left( \left\vert \Theta _{Z}\left( \mu \right) -\Gamma (\mu
)\right\vert >\varepsilon y;\mu <\Upsilon \right) \\
&\leq &\mathbf{E}\left[ \sum_{n=1}^{\infty }\mathbf{P}\left( \left\vert
\Theta _{Z}(\mu ,\mu +n)\mathbf{-}\left\langle \mathbf{Z}(\mu ,\mu +n),%
\mathbf{C}_{\mu +n}\right\rangle \right\vert >\frac{\varepsilon y}{2n^{2}}%
\Big|\, \mathcal{C}_{\mu }\right) I\left( \mu <\Upsilon \right) %
\right] \\
&\leq &\frac{m^{\kappa }K_{1}}{\varepsilon ^{\kappa }y^{\kappa }}\mathbf{E}%
\left[ \sum_{n=1}^{\infty }n^{2\kappa }\left\Vert \mathbf{V}(\mu
)\right\Vert ^{\kappa /2}\left\Vert \Pi _{\mu ,\mu +n}\right\Vert ^{\kappa
/2}I\left( \mu <\Upsilon \right) \right] \\
&=&\frac{m^{\kappa }K_{1}}{\varepsilon ^{\kappa }y^{\kappa }}\mathbf{E}\left[
\left\Vert \mathbf{V}(\mu )\right\Vert ^{\kappa /2}I\left( \mu <\Upsilon
\right) \sum_{n=1}^{\infty }n^{2\kappa }\mathbf{E}\left\Vert \Pi
_{0,n}\right\Vert ^{\kappa /2}\right] \\
&\leq &\frac{K_{2}}{\varepsilon ^{\kappa }y^{\kappa }r^{\kappa /2}}\mathbf{E}%
\left[ \left\Vert \mathbf{V}(\mu )\right\Vert ^{\kappa }I\left( \mu
<\Upsilon \right) \right] \leq \frac{\varepsilon }{y^{\kappa }}\mathbf{E}%
\left[ \left\Vert \mathbf{V}(\mu )\right\Vert ^{\kappa }I\left( \mu
<\Upsilon \right) \right]
\end{eqnarray*}%
for all $r\geq r(\varepsilon )$.

The lemma is proved.

Up to now we have assumed that the initial number of particles and the
initial amount of the final product are nonrandom. The next two lemmas are
free of this restriction.

Let%
\begin{equation*}
\Theta \left( N\right) :=\sum_{n=0}^{N}\Delta _{I}(n)
\end{equation*}

\begin{lemma}
\label{Lremaindd}Let a MBPIFPRE be subcritical, and $n=1$ be the
starting point of a life period of the process. If there exists
$\delta >0$ such that
\begin{equation}
\mathbf{E}\Theta ^{\kappa +\delta }(0)<\infty ,\quad \max_{1\leq i\leq m}%
\mathbf{E}\varphi _{i}^{\kappa +\delta }<\infty ,\qquad \mathbf{E}\psi
^{\kappa +\delta }<\infty ,  \label{finit}
\end{equation}%
then for any $\varepsilon >0$ there exists $r=r(\varepsilon )$ such that for
all $y\geq y_{0}$
\begin{equation}
\mathbf{P}\left( \Theta \left( \mu -1\right) >\varepsilon y;\mu <\Upsilon
\right) <\frac{1}{y^{\kappa +\delta /2}}.  \label{Dif33}
\end{equation}
\end{lemma}

\textbf{Proof}. Let, as earlier, $\Upsilon $ be the length of the
life period in question. Recalling (\ref{DDphi}) put%
\begin{equation*}
\Theta \left( N,r\right) :=\Theta (0)+\sum_{l=0}^{N}\left( \psi
(l+1)+\sum_{i=1}^{m}\sum_{k=1}^{r}\varphi _{i}(l;k)\right) .
\end{equation*}%
Clearly,%
\begin{eqnarray*}
&&\mathbf{P}\left( \Theta \left( N-1,r\right) >\varepsilon y;N<\Upsilon
\right) \\
&&\quad\leq \mathbf{P}\left( \Theta (0)>\frac{\varepsilon y}{3}\right) \mathbf{P}%
\left( N<\Upsilon \right) +\mathbf{P}\left( \sum_{l=0}^{N}\psi (l+1)>\frac{%
\varepsilon y}{3}\right) \\
&&\qquad\qquad+\mathbf{P}\left(
\sum_{l=0}^{N}\sum_{i=1}^{m}\sum_{k=1}^{r}\varphi
_{i}(l;k)>\frac{\varepsilon y}{3}\right) \\
&&\quad\leq\left( \frac{3}{\varepsilon y}\right) ^{\kappa +\delta }\mathbf{E}%
\Theta ^{\kappa +\delta }(0)\mathbf{P}\left( N<\Upsilon \right)
+\left(
N+1\right) \mathbf{P}\left( \psi >\frac{\varepsilon y}{3\left( N+1\right) }%
\right) \\
&&\qquad\qquad+\left( N+1\right) r\mathbf{P}\left( \sum_{i=1}^{m}\varphi _{i}>\frac{%
\varepsilon y}{3(N+1)r}\right) \\
&&\quad\leq\left( \frac{3}{\varepsilon y}\right) ^{\kappa +\delta }\mathbf{E}%
\Theta ^{\kappa +\delta }(0)\mathbf{P}\left( N<\Upsilon \right)
+\left( N+1\right) \left( \frac{3\left( N+1\right) }{\varepsilon
y}\right) ^{\kappa
+\delta }\mathbf{E}\psi ^{\kappa +\delta } \\
&&\qquad\qquad+\left( N+1\right) r\left( \frac{3\left( N+1\right) r}{\varepsilon y}%
\right) ^{\kappa +\delta }\mathbf{E}\left( \sum_{i=1}^{m}\varphi
_{i}\right) ^{\kappa +\delta }.
\end{eqnarray*}
In view of \ Lemma \ref{LtailAsympt} there exist constants $K_{3}$
and $K_{4}$  such that for $y>0$
\begin{equation*}
\sum_{N>K_{4}\ln y}\mathbf{P}\left( N<\Upsilon \right) \leq \frac{K_{3}}{%
y^{\kappa +\delta }}.
\end{equation*}%
Now
using conditions (\ref{finit}) we have%
\begin{eqnarray*}
\mathbf{P}\left( \Theta \left( \mu -1\right) >\varepsilon y;\mu <\Upsilon
\right) &=&\sum_{N=1}^{\infty }\mathbf{P}\left( \Theta \left( N-1\right)
>\varepsilon y;\mu =N<\Upsilon \right) \\
&\leq &\sum_{1\leq N\leq K_{4}\ln y}\mathbf{P}\left( \Theta \left(
N-1,r\right) >\varepsilon y;\mu =N<\Upsilon \right) \\
&&+\sum_{N>K_{4}\ln y}\mathbf{P}\left( N<\Upsilon \right) \\
&\leq &\frac{K_{5}}{\left( \varepsilon y\right) ^{\kappa +\delta }}%
\sum_{1\leq N\leq K_{4}\ln y}\left( N+1\right) ^{\kappa +\delta +1} \\
&\leq &\frac{K_{6}\ln ^{\kappa +\delta +2}y}{\left( \varepsilon y\right)
^{\kappa +\delta }}\leq \frac{1}{y^{\kappa +\delta /2}}.
\end{eqnarray*}%
The lemma is proved.

\begin{lemma}
\label{Lremtail}Under the conditions of Lemma \ref{Lremaindd} there exists $%
r $ such that for all $y\geq y_{0}$
\begin{equation}
\mathbf{P}\left( \Theta >y;\mu >\Upsilon \right) \leq \frac{1}{y^{\kappa
+\delta /2}}.  \label{Dif44}
\end{equation}
\end{lemma}

\textbf{Proof}. Letting $K_{3}$ and $K_{4}$ be the same as in the previous
lemma we have for $y\geq y_{0}$%
\begin{eqnarray*}
&&\mathbf{P}\left( \Theta >y;\mu >\Upsilon \right) \leq \mathbf{P}\left(
\Theta \left( \Upsilon ,r\right) >y\right) \\
&&\quad\leq\mathbf{P}\left( \Theta \left( \left[ K_{4}\ln y\right]
,r\right)
>y\right) +\mathbf{P}\left( \Upsilon >K_{4}\ln y\right) \\
&&\quad\leq \mathbf{P}\left( \Theta (0)>\frac{y}{3}\right)
+\sum_{l=0}^{\left[
K_{4}\ln y\right] }\mathbf{P}\left( \psi (l+1)>\frac{y}{3\left[ K_{4}\ln y%
\right] }\right) \\
&&\quad\quad+\sum_{l=0}^{\left[ K_{4}\ln y\right] }\sum_{i=1}^{m}\sum_{k=1}^{r}\mathbf{%
P}\left( \varphi _{i}(l;k)>\frac{y}{3mr\left[ K_{4}\ln y\right] }\right) +%
\mathbf{P}\left( \Upsilon >K_{4}\ln y\right) \\
&&\qquad\qquad\qquad\qquad\qquad\leq K\frac{\ln ^{\kappa +\delta +1}y}{y^{\kappa +\delta }}+\frac{1}{%
y^{\kappa +\delta }}.
\end{eqnarray*}%
The lemma is proved.

\section{Proof of Theorem \protect\ref{TTmain}}\label{MM}

Now we are ready to prove the main result of the paper, Theorem \ref{TTmain}%
. First observe that by the equivalence of the norms $\left\Vert \cdot
\right\Vert $ and $\left\Vert \cdot \right\Vert _{2}$ and estimates (\ref%
{Dif11}),(\ref{Dif22}), (\ref{Dif33}) and (\ref{Dif44}), for any $%
\varepsilon \in (0,1/3)$ one can find $r=r(\varepsilon )$ such that for all $%
y\geq y_{0}(r,\varepsilon )$
\begin{eqnarray}
&&\mathbf{P}\left( \Theta >y\right) \leq \mathbf{P}\left( \left\langle
\mathbf{V}(\mu ),\Xi _{\mu }\right\rangle >y(1-3\varepsilon );\mu <\Upsilon
\right)  \notag \\
&&  \notag \\
&&+\mathbf{P}\left( \left\vert \Gamma (\mu )-\left\langle \mathbf{V}(\mu
),\Xi _{\mu }\right\rangle \right\vert >\varepsilon y;\mu <\Upsilon \right) +%
\mathbf{P}\left( \left\vert \Theta _{Z}\left( \mu \right) -\Gamma (\mu
)\right\vert >\varepsilon y;\mu <\Upsilon \right)  \notag \\
&&  \notag \\
&&+\mathbf{P}\left( \Theta >y;\mu >\Upsilon \right) +\mathbf{P}\left( \Theta
\left( \mu -1\right) >\varepsilon y;\mu <\Upsilon \right)  \notag \\
&&  \notag \\
&&\leq \mathbf{P}\left( \left\langle \mathbf{V}(\mu ),\Xi _{\mu
}\right\rangle >y(1-3\varepsilon );\mu <\Upsilon \right)
+\frac{2\varepsilon }{y^{\kappa }}\mathbf{E}\left[ \left\Vert
\mathbf{V}(\mu )\right\Vert _{2}^{\kappa }I\left\{ \mu <\Upsilon
\right\} \right] +\frac{2}{y^{\kappa
+\delta /2}}.  \notag \\
&&  \label{FiDecompp}
\end{eqnarray}%
Let $l(\mathbf{u})>0$ be the function involved in Condition~$T$. By this
condition and the independency of $\mathbf{V}(\mu )$ and $\Xi _{\mu }$ given
$\mu <\Upsilon $ we conclude
\begin{eqnarray*}
&&\lim \sup_{y\rightarrow \infty }y^{\kappa }\mathbf{P}\left( \left\langle
\mathbf{V}(\mu ),\Xi _{\mu }\right\rangle >y(1-3\varepsilon );\mu <\Upsilon
\right) \\
&&\leq \lim \sup_{y\rightarrow \infty }y^{\kappa }\int_{\left\Vert \mathbf{u}%
\right\Vert _{2}=r}^{\infty }\mathbf{P}\left( \mathbf{V}(\mu )\in d\mathbf{u}%
;\mu <\Upsilon \right) \mathbf{P}\left( \left\langle \frac{\mathbf{u}}{%
\left\Vert \mathbf{u}\right\Vert _{2}},\Xi \right\rangle \geq \frac{%
(1-3\varepsilon )y}{\left\Vert \mathbf{u}\right\Vert _{2}}\right) \\
&&=K_{0}(1-3\varepsilon )^{-\kappa }\int_{\left\Vert
\mathbf{u}\right\Vert _{2}=r}^{\infty }\mathbf{P}\left(
\mathbf{V}(\mu )\in d\mathbf{u};\mu <\Upsilon \right) \left\Vert
\mathbf{u}\right\Vert _{2}^{\kappa }l\left(
\frac{\mathbf{u}}{\left\Vert \mathbf{u}\right\Vert _{2}}\right) \\
&&=K_{0}(1-3\varepsilon )^{-\kappa }\mathbf{E}\left[ \left\Vert \mathbf{V}%
(\mu )\right\Vert _{2}^{\kappa }l\left( \frac{\mathbf{V}(\mu
)}{\left\Vert
\mathbf{V}(\mu )\right\Vert _{2}}\right) I\left\{ \mu <\Upsilon \right\} %
\right] <\infty .
\end{eqnarray*}%
Since $y^{\kappa }\mathbf{P}\left( \Theta >y\right) $ does not depend on $r$
and $\varepsilon ,$ the previous estimate and (\ref{FiDecompp}) yield
\begin{equation}
\lim \sup_{y\rightarrow \infty }y^{\kappa }\mathbf{P}\left( \Theta >y\right)
<\infty  \label{PPinff}
\end{equation}%
and, moreover,%
\begin{equation}
\lim \sup_{y\rightarrow \infty }y^{\kappa }\mathbf{P}\left( \Theta >y\right)
\leq K_{0}\lim \inf_{r\rightarrow \infty }\mathbf{E}\left[ \left\Vert
\mathbf{V}(\mu )\right\Vert _{2}^{\kappa }l\left( \frac{\mathbf{V}(\mu )}{%
\left\Vert \mathbf{V}(\mu )\right\Vert _{2}}\right) I\left\{ \mu <\Upsilon
\right\} \right] .  \label{PP22}
\end{equation}%
To get a similar estimate from below we use for $\varepsilon >0$ the
inequality%
\begin{eqnarray*}
\mathbf{P}\left( \Theta >y\right) &\geq &\mathbf{P}\left( \Theta _{Z}\left(
\mu \right) ;\mu <\Upsilon \right) \\
&\geq &\mathbf{P}\left( \left\langle \mathbf{V}(\mu ),\Xi _{\mu
}\right\rangle >y(1+2\varepsilon );\mu <\Upsilon \right) -\mathbf{P}\left(
\left\vert \Theta _{Z}\left( \mu \right) -\Gamma (\mu )\right\vert
>\varepsilon y;\mu <\Upsilon \right) \\
&&-\mathbf{P}\left( \left\vert \Gamma \left( \mu \right) -\left\langle
\mathbf{V}(\mu ),\Xi _{\mu }\right\rangle \right\vert >\varepsilon y;\mu
<\Upsilon \right)
\end{eqnarray*}%
Now we select $r$ as large to meet estimates (\ref{Dif11}), (\ref{Dif22})
and (\ref{Dif33}). This gives for sufficiently large $y>r$ the inequality%
\begin{equation*}
\mathbf{P}\left( \Theta >y\right) \geq \mathbf{P}\left( \left\langle \mathbf{%
V}(\mu ),\Xi _{\mu }\right\rangle >y(1+2\varepsilon );\mu <\Upsilon \right) -%
\frac{2\varepsilon }{y^{\kappa }}\mathbf{E}\left[ \left\Vert \mathbf{V}(\mu
)\right\Vert _{2}^{\kappa }I\left\{ \mu <\Upsilon \right\} \right] .
\end{equation*}%
Letting $y$ tend to infinity we obtain%
\begin{eqnarray*}
&&\lim \inf_{y\rightarrow \infty }y^{\kappa }\mathbf{P}\left( \left\langle
\mathbf{V}(\mu ),\Xi _{\mu }\right\rangle >y(1+2\varepsilon );\mu <\Upsilon
\right) \\
&&=\lim \inf_{y\rightarrow \infty }y^{\kappa }\int_{\left\Vert \mathbf{u}%
\right\Vert _{2}=r}^{\infty }\mathbf{P}\left( \mathbf{V}(\mu )\in d\mathbf{u}%
;\mu <\Upsilon \right) \mathbf{P}\left( \left\langle \frac{\mathbf{u}}{%
\left\Vert \mathbf{u}\right\Vert _{2}},\Xi \right\rangle \geq \frac{%
(1+2\varepsilon )y}{\left\Vert \mathbf{u}\right\Vert _{2}}\right) \\
&&\qquad\qquad\qquad\qquad\qquad=K_{0}(1+2\varepsilon )^{-\kappa }\mathbf{E}\left[ \left\Vert \mathbf{V}%
(\mu )\right\Vert _{2}^{\kappa }l\left( \frac{\mathbf{V}(\mu
)}{\left\Vert
\mathbf{V}(\mu )\right\Vert _{2}}\right) I\left\{ \mu <\Upsilon \right\} %
\right] .
\end{eqnarray*}%
Let $K(l):=\inf_{\mathbf{u\in U}_{+}}l(\mathbf{u})>0$. We know that for
sufficiently small $\varepsilon >0$ and an appropriate $r$%
\begin{eqnarray*}
&&K_{0}(1+2\varepsilon )^{-\kappa }\mathbf{E}\left[ \left\Vert \mathbf{V}%
(\mu )\right\Vert _{2}^{\kappa }l\left( \frac{\mathbf{V}(\mu )}{\left\Vert
\mathbf{V}(\mu )\right\Vert _{2}}\right) I\left\{ \mu <\Upsilon \right\} %
\right] -2\varepsilon \mathbf{E}\left[ \left\Vert \mathbf{V}(\mu
)\right\Vert _{2}^{\kappa }I\left\{ \mu <\Upsilon \right\} \right] \\
&&\qquad\qquad\qquad\qquad\qquad\geq \left(
K_{0}K(l)(1+2\varepsilon )^{-\kappa }-2\varepsilon \right)
\mathbf{E}\left[ \left\Vert \mathbf{V}(\mu )\right\Vert
_{2}^{\kappa }I\left\{ \mu <\Upsilon \right\} \right] >0
\end{eqnarray*}%
which implies
\begin{equation*}
\lim \inf_{y\rightarrow \infty }y^{\kappa }\mathbf{P}\left( \Theta >y\right)
>0
\end{equation*}%
leading in turn to%
\begin{equation*}
\lim \inf_{y\rightarrow \infty }y^{\kappa }\mathbf{P}\left( \Theta >y\right)
\geq K_{0}\lim \sup_{r\rightarrow \infty }\mathbf{E}\left[ \left\Vert
\mathbf{V}(\mu )\right\Vert _{2}^{\kappa }l\left( \frac{\mathbf{V}(\mu )}{%
\left\Vert \mathbf{V}(\mu )\right\Vert _{2}}\right) I\left\{ \mu <\Upsilon
\right\} \right] .
\end{equation*}%
This combined with (\ref{PP22}) shows existence of the limit%
\begin{equation*}
\lim_{r\rightarrow \infty }\mathbf{E}\left[ \left\Vert \mathbf{V}(\mu
)\right\Vert _{2}^{\kappa }l\left( \frac{\mathbf{V}(\mu )}{\left\Vert
\mathbf{V}(\mu )\right\Vert _{2}}\right) I\left\{ \mu <\Upsilon \right\} %
\right] \in (0,\infty )
\end{equation*}%
and gives
\begin{equation*}
\lim_{y\rightarrow \infty }y^{\kappa }\mathbf{P}\left( \Theta >y\right)
=K_{0}\lim_{r\rightarrow \infty }\mathbf{E}\left[ \left\Vert \mathbf{V}(\mu
)\right\Vert _{2}^{\kappa }l\left( \frac{\mathbf{V}(\mu )}{\left\Vert
\mathbf{V}(\mu )\right\Vert _{2}}\right) I\left\{ \mu <\Upsilon \right\} %
\right]
\end{equation*}
proving (\ref{EstTail}). The validity of (\ref{EstMom}) is now
evident.

The theorem is proved.

\section{Polling systems\label{Spoll}}

We consider a polling system consisting of a single server and $m$
stations with infinite-buffer queues indexed by $i\in \left\{
1,\ldots ,m\right\} .$ Initially there are no customers in the
system and the server waits their arrival at a parking place $R$.
Customers arrive to the queues in accordance with a point process
whose parameters are changing in a random manner each time when
the server switches from station to station. When a customer
arrives to the system, say to station $i$, the server  immediately
begins the service by visiting the stations in cyclic order
$(i\rightarrow i+1\rightarrow \cdots \rightarrow m\rightarrow
1\rightarrow \cdots )$ starting at station $i$ according to a
selected service policy (to be described later on) and with random
positive switch-over times between the queues. Later on the
initial stage of services $(i\rightarrow i+1\rightarrow \cdots
\rightarrow m\rightarrow \cdots \rightarrow i-1)$ will be called
the zero cycle. The subsequent routes of the server will be called
the first cycle, the second cycle and so on. \ The server
immediately moves to its parking place if there are no customers
in the system to be served.

Our goal is to investigate various characteristics related with busy periods
of the system whose work starts by the arrival of a single customer to a
station $i\in \left\{ 1,...,m\right\} $ at moment 0$.$

To give a rigorous description of the arrival and service processes for the
system in question we need some notions. Let
\begin{equation*}
\chi ^{(i)}(\mathbf{s;}\lambda ):=\mathbf{E}\left[ s_{1}^{\theta
_{i1}}s_{2}^{\theta _{i2}}\cdots s_{m}^{\theta _{im}}e^{-\lambda \phi _{i}}%
\right], i=1,...,m
\end{equation*}
be the mixed probability generating functions (m.p.g.f.'s) of
$(m+1)-$ dimensional vectors $\left( \theta _{i1},\ldots ,\theta
_{im};\phi _{i}\right) $, where $\theta _{ij}$ are nonnegative
integer-valued random
variables and $\phi _{i}$ is a nonnegative random variable. Denote%
\begin{equation*}
\mathbf{\chi }(\mathbf{s;}\lambda ):=\left( \chi
^{(1)}(\mathbf{s;}\lambda ),\ldots ,\chi ^{(m)}(\mathbf{s;}\lambda
)\right)
\end{equation*}%
the respective vector-valued m.p.g.f.

Let, further,
\begin{equation*}
\rho ^{(i)}(\mathbf{s;}\lambda ):=\mathbf{E}\left[ s_{1}^{\zeta
_{i1}}s_{2}^{\zeta _{i2}}\cdots s_{m}^{\zeta _{im}}e^{-\lambda \sigma _{i}}%
\right] ,\,\ i=1,...,m
\end{equation*}%
be the m.p.g.f.'s of $(m+1)-$ dimensional vectors $\left( \zeta _{i1},\ldots
,\zeta _{im};\sigma _{i}\right) $, where $\zeta _{ij}$ are nonnegative
integer-valued random variables and $\sigma _{i}$ is a nonnegative random
variable. Denote%
\begin{equation*}
\mathbf{\rho }(\mathbf{s;}\lambda ):=\left( \rho ^{(1)}(\mathbf{s;}\lambda
),\ldots ,\rho ^{(m)}(\mathbf{s;}\lambda )\right)
\end{equation*}%
the respective vector-valued m.p.g.f.\ and let $\mathcal{H}=\left\{ (\mathbf{%
\chi }(\mathbf{s;}\lambda ),\mathbf{\rho }(\mathbf{s;}\lambda
))\right\} $ be the set of pairs of the described m.p.g.f.'s.
Assume that a probability
measure $\mathbb{P}$ is specified on the natural $\sigma $-algebra of $%
\mathcal{H}$ and let
\begin{equation*}
(\mathbf{\chi }_{0}(\mathbf{s;}\lambda ),\mathbf{\rho }_{0}(\mathbf{s;}%
\lambda )),(\mathbf{\chi }_{1}(\mathbf{s;}\lambda ),\mathbf{\rho }_{1}(%
\mathbf{s;}\lambda )),(\mathbf{\chi }_{2}(\mathbf{s;}\lambda ),\mathbf{\rho }%
_{2}(\mathbf{s;}\lambda )),\ldots ,
\end{equation*}%
be a sequence of pairs $(\mathbf{\chi }_{n}(\mathbf{s;}\lambda ),\mathbf{%
\rho }_{n}(\mathbf{s;}\lambda ))$ selected in iid manner from $\mathcal{H}$
according to measure $\mathbb{P},$ where $\mathbf{\chi }_{n}(\mathbf{s;}%
\lambda ):=\left( \chi _{n}^{(1)}(\mathbf{s;}\lambda ),\ldots
,\chi
_{n}^{(m)}(\mathbf{s;}\lambda )\right) $ and $\mathbf{\rho }_{n}(\mathbf{s;}\lambda ):=\left( \rho _{n}^{(1)}(\mathbf{%
s;}\lambda ),\ldots ,\rho _{n}^{(m)}(\mathbf{s;}\lambda )\right) $
have the components
\begin{equation*}
\chi _{n}^{(i)}(\mathbf{s;}\lambda ):=\mathbf{E}\left[
s_{1}^{\theta _{i1}(n)}s_{2}^{\theta _{i2}(n)}\cdots s_{m}^{\theta
_{im}(n)}e^{-\lambda \phi _{i}(n)}\right],
\end{equation*}%
\begin{equation*}
\rho _{n}^{(i)}(\mathbf{s;}\lambda ):=\mathbf{E}\left[ s_{1}^{\zeta
_{i1}(n)}s_{2}^{\zeta _{i2}(n)}\cdots s_{m}^{\zeta _{im}(n)}e^{-\lambda
\sigma _{i}(n)}\right] .
\end{equation*}

In the present paper we investigate such polling systems whose arrival and
service procedures of customers meet the following two conditions.

\textbf{Branching property}. When the server arrives to station $i$ for the $%
n$-th time and find there, say, $k_{i}$ customers labelled $1,2,...,k_{i}$ ,
then, during the course of the server's visit, the arrival of customers is
arranged in such a way that after the end of each stage of service of
customer $j$ (the number of such stages may be more than one if the service
discipline admits feedback) the queues at the system will be increased by a
random population of customers $\left( \theta _{i1}(n,j),\ldots ,\theta
_{im}(n,j)\right) ,$ where $\theta _{il}(n,j)$ is the number of customers
added to the $l$-th station, and, in addition, a final product of size $\phi
_{i}(n,j)\geq 0$ will be added to the system. It is assumed that the vectors
$\left( \theta _{i1}(n,j),\ldots ,\theta _{im}(n,j);\phi _{i}(n,j)\right) ,$
$j=1,2,\ldots ,k_{i}$ are iid and such that%
\begin{equation*}
\mathbf{E}\left[ s_{1}^{\theta _{i1}(n,j)}s_{2}^{\theta _{i2}(n,j)}\cdots
s_{m}^{\theta _{im}(n,j)}e^{-\lambda \phi _{i}(n,j)}\right] =\chi _{n}^{(i)}(%
\mathbf{s;}\lambda ).
\end{equation*}

Note that $\chi _{n}^{(i)}(\mathbf{s;}\lambda )$ is a random
m.p.g.f. and, therefore, the parameters of the polling system are
changed from cycle to cycle in a random manner.

\textbf{Immigration property}. Upon departing from the $i-$th station for the $%
n-$th time the server switches to the queue existing at station
$(i+1)(mod\,m)$. The process of switching takes a random time at
the end of which the queues at the system will be increased by a
random population of customers $\left( \zeta _{i1}(n),\ldots
,\zeta _{im}(n)\right) ,$ where $\zeta _{il}(n)$ is the number of
customers added to the $l$-th station, and, in addition, a final
product of size $\sigma _{i}(n)\geq 0$ will be added to the
system. It is assumed that the vectors $\left( \zeta
_{i1}(n),\ldots ,\zeta
_{im}(n);\sigma _{i}(n)\right) ,n=1,2,...$ are iid and such that%
\begin{equation*}
\mathbf{E}\left[ s_{1}^{\zeta _{i1}(n)}s_{2}^{\zeta _{i2}(n)}\cdots
s_{m}^{\zeta _{im}(n)}e^{-\lambda \sigma _{i}(n)}\right] =\rho _{n}^{(i)}(%
\mathbf{s;}\lambda ).
\end{equation*}

In the sequel we call the described polling systems as the branching type
polling systems with final product and random environment (BTPSFPRE).

Polling systems are rather popular subject of investigations in
queueing theory (see surveys \cite{ViSe06} , \cite{Mei07} and
\cite{YY07} for definitions and more details). For instance,
polling systems possessing
the branching property in which the probability generating functions $%
h_{n}^{(i)}(\mathbf{s}),\quad n=1,2,\ldots $ are nonrandom and
identical for all cycles were considered in \cite{Fur91},
\cite{Res93}, \cite{VD2002} and quite recently in \cite{Alt2009}
and \cite{AltF2007}. These models cover many classical service
policies, including the exhaustive, gated, binomial-gated and
their feedback modifications.

Polling systems with input parameters and service disciplines changing in a
random manner and (or) depending on the states of systems were investigated
by the fluid method in \cite{FC98}-\cite{FL98}) and by the method based on
the construction of appropriate Lyapunov functions in \cite{PMPV2006} - \cite%
{PMPP2008}. The present paper may be viewed as a compliment to article \cite%
{Vat2009} in which the BTPSFPRE with zero switchover times were
analyzed by means of MBPRE without immigration.

Before we pass to general results, consider two examples of BTPSFPRE.

Let $\mathcal{T}_{+}\mathcal{=}\left\{ T\right\} $ be the set of all
probability distributions of nonnegative random variables,
\begin{equation*}
\mathcal{T}_{+}^{m}\times \mathcal{I}_{+}^{m}:=\left\{ \left(
T_{1},...,T_{m}\right) \times \left( I_{1},...,I_{m}\right) :T_{i}\in
\mathcal{T}_{+},I_{j}\in \mathcal{I}_{+},i,j=1,...,m\right\}
\end{equation*}%
be the set of all $2m$-dimensional tuples of such distributions, $\mathcal{M}%
_{\varepsilon }=\left\{ \mathcal{E}\right\} $ be the set of all $m\times m$
matrices $\mathcal{E}=\left( \varepsilon _{ij}\right) _{i,j=1}^{m}$ with
nonnegative elements, and $\mathcal{M}_{\gamma }=\left\{ \Gamma \right\} $
be the set of all $m\times (m+1)$ matrices $\Gamma =\left( \gamma
_{ij}\right) _{i=1,j=0}^{m}$ with nonnegative elements such that

\begin{equation*}
\sum_{j=0}^{m}\gamma _{ij}=1,\quad i=1,...,m,\quad \max_{1\leq i\leq
m}\gamma _{i0}>0.
\end{equation*}

Let $\mathbf{P}$\ be a measure on the Borel $\sigma -$algebra of the space
\begin{equation*}
\mathcal{S}:=\mathcal{M}_{\varepsilon }\times \mathcal{M}_{\varepsilon
}\times \mathcal{M}_{\gamma }\times \mathcal{T}_{+}^{m}\times \mathcal{I}%
_{+}^{m}.
\end{equation*}

\begin{example}
\label{Examp1}(motivated by \cite{PMPP2008}). Consider a polling system with
$m$ stations and a single server performing cyclic service of the customers
at the stations. Assume that initially there are no customers in the system
and the server is located at a parking place $R$. Assume that given the idle
system the flow of customers arriving to station $l$ is Poisson with, say, a
deterministic rate $\varepsilon _{l}$. When the first customer appears in
the system the server selects a random element $\left( \mathcal{E}_{0},%
\mathcal{E}_{0}^{J},\Gamma _{0},\mathbf{T}_{0}\times \mathbf{I}_{0}\right)
\in \mathcal{S}$ with%
\begin{eqnarray*}
\mathcal{E}_{0} &\mathcal{=}&\left( \varepsilon _{ij}(0)\right)
_{i,j=1}^{m},\, \mathcal{E}_{0}^{J}\mathcal{=}\left( \varepsilon
_{ij}^{I}(0)\right) _{i,j=1}^{m},\quad \Gamma _{0}=\left( \gamma
_{ij}(0)\right) _{i=1,j=0}^{m},
\\
\quad \mathbf{T}_{0} &=&\left( T_{10},...,T_{m0}\right) ,\quad \,\mathbf{I}%
_{0}=\left( I_{10},...,I_{m0}\right)
\end{eqnarray*}%
and immediately starts its zero service cycle $(i\rightarrow i+1(mod
m)\rightarrow \cdots \rightarrow i-1(mod m))$ adopting the gated server
policy with random switch-over times. Namely, the server serves all the
customers that were queueing, say at station $i$ when the server arrived and
then proceeds to the next (in cyclic order) station with switch-over times $%
\sigma _{i}(0)$ distributed according to
$I_{i0}(x):=\mathbf{P}\left( \sigma _{i}(0)\leq x\right) $.
During the switch-over period $\sigma _{i}(0)$ new customers
arrive to the system
according to independent Poisson flows with intensities given by the vector $%
(\varepsilon _{i1}^{I}(0),\varepsilon _{i2}^{I}(0),...,\varepsilon
_{im}^{I}(0))$ (some of the components may be equal to zero).

 For
the period while the server performs the batch of services at
station $i$ new customers arrive to the system according to
independent Poisson flows with intensities given by the vector
$(\varepsilon _{i1}(0),\varepsilon _{i2}(0),...,\varepsilon
_{im}(0))$ (some of the components may be equal to zero) and the
service times of customers are iid and distributed according to
$T_{i0}(x):=\mathbf{P}\left( \tau _{i}(0)\leq x\right) $. Each
served customer either goes to station $j\in \{1,\ldots ,m\} $
with probability $\gamma _{ij}(0)$ or leaves the system with
probability $\gamma _{i0}(0)$ independently of other events.
Besides, after the end of each service stage of a customer the
customer contributes\ to the system the service time at this stage
as the final product.

It is assumed that given $\left( \mathcal{E}_{0},\mathcal{E}_{0}^{J},\Gamma
_{0},\mathbf{T}_{0}\times \mathbf{I}_{0}\right) $ the service times,
switch-over times and the arrival process of new customers are independent.

The subsequent routes $n=1,2,\ldots $ have the same probabilistic structure
specified by the tuples
\begin{eqnarray*}
\mathcal{E}_{n} &\mathcal{=}&\left( \varepsilon _{ij}(n)\right)
_{i,j=1}^{m}, \,
 \mathcal{E}_{n}^{J}\mathcal{=}\left( \varepsilon
_{ij}^{I}(n)\right) _{i,j=1}^{m}, \quad\Gamma _{n}=\left( \gamma
_{ij}(n)\right) _{i=1,j=0}^{m},
\\
\quad \mathbf{T}_{n} &=&\left( T_{1n},...,T_{mn}\right) ,\,\quad \mathbf{I}%
_{n}=\left( I_{1n},...,I_{mn}\right)
\end{eqnarray*}%
with only difference that at the beginning of cycle $n\geq 1$ there is a
possibility to have more than one customer in the system.
\end{example}

One can show (see \cite{Vat2009}) that this system possesses a branching property in which the m.p.g.f. $\mathbf{%
\chi }_{n}(\mathbf{s;}\lambda )$ has the components (in our
setting and the service time of customers as the final product)
\begin{eqnarray}
&&\chi _{n}^{(i)}(\mathbf{s;}\lambda ) :=\mathbf{E}\left[
s_{1}^{\theta _{i1}(n)}s_{2}^{\theta _{i2}(n)}\cdots s_{m}^{\theta
_{im}(n)}e^{-\lambda
\tau _{i}(n)}\right]  \notag \\
&&\qquad\qquad\,=\left( \gamma _{i0}(n)+\sum_{j=1}^{m}\gamma
_{ij}(n)s_{j}\right) t_{in}\left(\lambda
+\sum_{j=1}^{m}\varepsilon _{ij}(n)(1-s_{j})\right),
\label{GateTotal}
\end{eqnarray}%
where%
\begin{equation*}
t_{in}(\lambda ):=\int_{0}^{\infty }e^{-\lambda x}dT_{in}(x),
\end{equation*}%
and an immigration property with $\mathbf{\rho }_{n}(\mathbf{s;}\lambda )$
having components%
\begin{eqnarray}
&&\rho _{n}^{(i)}(\mathbf{s;}\lambda ) :=\mathbf{E}\left[
s_{1}^{\zeta _{i1}(n)}s_{2}^{\zeta _{i2}(n)}\cdots s_{m}^{\zeta
_{im}(n)}e^{-\lambda
\sigma _{i}(n)}\right]  \notag \\
&&\qquad\qquad\,=\int_{0}^{\infty }\exp \left\{ -\left( \lambda
+\sum_{j=1}^{m}\varepsilon _{ij}^{I}(n)(1-s_{j})\right) x\right\}
dI_{in}(x).  \label{ImmigrTotal}
\end{eqnarray}

\begin{example}
\label{Examp2}(compare with \cite{PMPP2008}). Consider the same
polling system as earlier but assume now that at each station the
server adopts the exhaustive server policy: it serves all the
customers that were queueing at the station when the server
arrived together with all subsequent arrivals up until the queue
at this station becomes empty and then switches over to the next
station.
\end{example}

Let%
\begin{equation}
w_{in}(\lambda ):=\frac{1-\gamma _{ii}(n)}{1-\gamma _{ii}(n)t_{in}(\lambda )}%
t_{in}(\lambda )  \label{ForW}
\end{equation}%
and%
\begin{equation*}
y_{in}(\mathbf{s}):=\frac{\gamma _{i0}(n)+\sum_{j\neq i}\gamma _{ij}(n)s_{j}%
}{1-\gamma _{ii}(n)}.
\end{equation*}

It is know (see, for instance, \cite{Vat2009}) that the system above
possesses a branching property with $\mathbf{\phi }_{n}(\mathbf{s;}\lambda )$
whose components are unique solutions of the equations%
\begin{equation*}
\chi _{n}^{(i)}(\mathbf{s;}\lambda
)=y_{in}(\mathbf{s})w_{in}\left( \lambda +\sum_{j\neq
i}\varepsilon _{ij}(n)(1-s_{j})+\varepsilon _{ii}(n)(1-\chi
_{n}^{(i)}(\mathbf{s;}\lambda ))\right) .
\end{equation*}%
The immigration property for this system fulfills with $\mathbf{\rho }_{n}(%
\mathbf{s;}\lambda )$ specified by (\ref{ImmigrTotal}).

Now we come back to the general case. We assume that the server
may start the service of a customer of any type arrived to an idle
system immediately and would like to study the distribution of the
total size of the final product accumulated in the system during a
busy period of the server. For this reason the law of arrival of
customers to an idle system plays no role for the subsequent
arguments. The only assumption we need that the probability of
arrival two or more customer to an idle system simultaneously is
zero.

Set%
\begin{equation*}
h_{n}^{(i)}(\mathbf{s}):=\chi
_{n}^{(i)}(\mathbf{s};0),i=1,2,\ldots ,m
\end{equation*}%
and for $n=0,1,2,\ldots $ introduce m.p.g.f.'s
\begin{equation*}
F_{n}^{(i)}(\mathbf{s;}\lambda )=\mathbf{E}\left[ s_{1}^{\xi
_{i1}(n)}s_{2}^{\xi _{i2}(n)}\cdots s_{m}^{\xi _{im}(n)}e^{-\lambda \varphi
_{i}(n)}\right]
\end{equation*}%
and p.g.f.'s
\begin{equation*}
f_{n}^{(i)}(\mathbf{s})=\mathbf{E}\left[ s_{1}^{\xi _{i1}(n)}s_{2}^{\xi
_{i2}(n)}\cdots s_{m}^{\xi _{im}(n)}\right]
\end{equation*}%
by the equalities $F_{n}^{(m)}(\mathbf{s};\lambda ):=\chi _{n}^{(m)}(\mathbf{%
s;}\lambda ),$
\begin{equation}
F_{n}^{(i)}(\mathbf{s};\lambda ):=\chi _{n}^{(i)}\left(
s_{1},\ldots ,s_{i},F_{n}^{(i+1)}\left( \mathbf{s};\lambda \right)
,\ldots ,F_{n}^{(m)}\left( \mathbf{s};\lambda \right) ;\lambda
\right) ,\,i<m, \label{GG3}
\end{equation}%
and $f_{n}^{(m)}(\mathbf{s}):=h_{n}^{(m)}(\mathbf{s}),$
\begin{equation}
f_{n}^{(i)}(\mathbf{s}):=h_{n}^{(i)}\left( s_{1},\ldots
,s_{i},f_{n}^{(i+1)}\left( \mathbf{s}\right) ,\ldots ,f_{n}^{(m)}\left(
\mathbf{s}\right) \right) ,\,i<m.  \label{GGGG}
\end{equation}%
Further, set
\begin{equation*}
l_{n}^{(i)}(\mathbf{s}):=\rho _{n}^{(i)}(\mathbf{s};0),i=1,2,\ldots ,m
\end{equation*}%
and for $n=0,1,2,\ldots $ introduce m.p.g.f.'s
\begin{eqnarray}
&&G_{n}(\mathbf{s;}\lambda ) :=\mathbf{E}\left[ s_{1}^{\eta
_{1}(n)}s_{2}^{\eta _{2}(n)}\cdots s_{m}^{\eta _{m}(n)}e^{-\lambda \psi (n)}%
\right]  \notag \\
&&\qquad\qquad=\prod_{i=1}^{m}\rho _{n}^{(i)}\left( s_{1},\ldots
,s_{i},F_{n}^{(i+1)}\left( \mathbf{s};\lambda \right) ,\ldots
,F_{n}^{(m)}\left( \mathbf{s};\lambda \right) ;\lambda \right)
\label{II2}
\end{eqnarray}%
and p.g.f.'s
\begin{eqnarray}
&&g_{n}(\mathbf{s}) :=\mathbf{E}\left[ s_{1}^{\eta
_{1}(n)}s_{2}^{\eta
_{2}(n)}\cdots s_{m}^{\eta _{m}(n)}\right]  \label{II3} \\
&&\qquad\quad=\prod_{i=1}^{m}l_{n}^{(i)}\left( s_{1},\ldots
,s_{i},f_{n}^{(i+1)}\left( \mathbf{s}\right) ,\ldots
,f_{n}^{(m)}\left( \mathbf{s}\right) \right) . \notag
\end{eqnarray}

We would like to describe conditions on the branching type polling system
under which power moments of the amount of the final product accumulated in
the system during its busy period are finite or infinite. Note that letting
the final product $\phi _{i}(n,j)$ be the service time of the $j$-th
customer served during the $n$-th visit of the server to station $i$ and $%
\sigma _{i}(n)$ be the switch-over time from station $i$ to station $%
(i+1)\,(mod\,m)$ after the $n$-th visit of the server to station $i$ we
provide conditions under which the tail distribution of the length of the
busy period of a polling system decays, as $y\rightarrow \infty ,$ like $%
y^{-\kappa }$ for some $\kappa >0$.

Our results are based on an important statement revealing
connections between the behavior of certain characteristics of the
busy periods of BTPSFPRE and the related characteristics of
life-length periods of a MBPIFPRE. To formulate this statement we
introduce the notion of generalized busy period as follows.

If a busy period of a BTPSFPRE starts by the arrival of a single
customer at station $J$ then, after the end of the busy period
generated by this customer we  continue to follow  the performance
of the system (with the described laws of arrival of customers
within switch-over times, the service disciplines and accumulating
of final product) up to the first visit of the server to station
$J$ when there are no customers in the system.

\begin{definition}
The time interval between the start of service of the first
customer and the moment when,  for the first time, the system is
idle after the end of the switch-over time from station $J-1$ to
station  $J$ is called a generalized busy period.
\end{definition}
 Denote by $\mathcal{M}_{J}$ the length of the respective
generalized busy period.

Since a generalized busy period includes a possibility to have an
idle system at the moment when the server arrives to station
$J-1(mod\,m)$ but have customers at the moment when the server
arrives to station $J$ (which arrive to the system during the
switch-over time $J-1\rightarrow J$), the length of a generalized
busy period is not less than the length of the standard busy
period.

In the sequel we assume without loss of generality that $J=1$.

\begin{theorem}
\label{RResRep}If a generalized busy period of a BTPSFPRE starts
by the arrival of a single customer at station $1$ then the
distribution of the amount of the final product accumulated in the
system to the end of the generalized period coincides with the
distribution of the total amount of the final product produced in
a MBPIFPRE within a life period which starts at moment $0$ by the
birth of a single particle of type $1$ and where the joint
distribution of the number of direct descendants and immigrants
and the amount of the final product produced by particles of
different types of the $k$-th generation is given by the
vector-valued m.p.g.f.'s
\begin{equation*}
\text{ }\mathbf{H}_{k}(\mathbf{s;}\lambda ):=\left( F_{k}^{(1)}(\mathbf{s;}%
\lambda ),\ldots ,F_{k}^{(m)}(\mathbf{s;}\lambda );G_{k+1}(\mathbf{s;}%
\lambda )\right) ,\,k=1,2,\ldots ,n
\end{equation*}%
specified by (\ref{GG3}) and (\ref{II2}).
\end{theorem}

To prove this theorem one should modify in a natural and evident way the
proof of Theorem 3 in \cite{Res93} and we omit the respective arguments.

The MBPIFPRE described in Theorem \ref{RResRep} will be called the \textit{%
associated}  for the BTPSFPRE.

Now we apply the obtained earlier  results for MBPIFPRE to our
polling system.

Let $H_{n}:=\left( h_{ij}(n)\right) _{i,j=1}^{m}$ be the matrix with elements%
\begin{equation*}
h_{ij}(n):=\frac{\partial h_{n}^{(i)}(\mathbf{s})}{\partial
s_{j}}\Big| _{\mathbf{s}=\mathbf{1}}
=\mathbf{E}_{\mathbf{h}_{n}}\theta _{ij}(n),
\end{equation*}%
and $A_{n}:=\left( a_{ij}(n)\right) _{i,j=1}^{m}$ be the matrix with elements%
\begin{equation}
a_{ij}(n):=\frac{\partial f_{n}^{(i)}(\mathbf{s})}{\partial
s_{j}}\Big| _{\mathbf{s}=\mathbf{1}}
=\mathbf{E}_{\mathbf{f}_{n}}\xi _{ij}(n). \label{DefElem}
\end{equation}%
Then in view of (\ref{GGGG}) $a_{mj}(n)=h_{mj}(n),j=1,2,\ldots ,m$ and for $%
i<m$
\begin{equation}
a_{ij}(n)=h_{ij}(n)I\left\{ j\leq i\right\}
+\sum_{k=i+1}^{m}h_{ik}(n)a_{kj}(n).  \label{lMatem}
\end{equation}%
For $i=1,\ldots ,m$ introduce auxiliary matrices%
\begin{equation*}
H_{n}^{(i)}:=\left(
\begin{array}{ccccccccc}
1 & 0 & \cdots & 0 & 0 & 0 & \cdots & \cdots & 0 \\
0 & 1 & \overset{\cdot }{}\cdot \underset{\cdot }{} & \overset{\cdot }{%
\underset{\cdot }{\cdot }} & 0 & 0 & \cdots & \cdots & 0 \\
\overset{\cdot }{\underset{\cdot }{\cdot }} & \overset{\cdot }{}\cdot
\underset{\cdot }{} & \overset{\cdot }{}\cdot \underset{\cdot }{} & 0 & 0 & 0
& \cdots & \cdots & 0 \\
0 & \cdots & 0 & 1 & 0 & 0 & \cdots & \cdots & 0 \\
h_{i1}(n) & h_{i2}(n) & \cdots & h_{i(i-1)}(n) & h_{ii}(n) & h_{i(i+1)}(n) &
\overset{\cdot }{\underset{\cdot }{\cdot }} & \overset{\cdot }{\underset{%
\cdot }{\cdot }} & h_{im}(n) \\
0 & \cdots & \cdots & 0 & 0 & 1 & 0 & \cdots & 0 \\
0 & \cdots & \cdots & 0 & 0 & 0 & 1 & \overset{\cdot }{}\cdot \underset{%
\cdot }{} & 0 \\
0 & \cdots & \cdots & 0 & 0 & 0 & \overset{\cdot }{}\cdot \underset{\cdot }{}
& \overset{\cdot }{}\cdot \underset{\cdot }{} & 0 \\
0 & \cdots & \cdots & 0 & 0 & 0 & \cdots & 0 & 1%
\end{array}%
\right) ,
\end{equation*}%
where for each $i$ the elements of the matrix $H_{n}$ are located only in
row $i$ of the matrix $H_{n}^{(i)}$. It is not difficult to check by (\ref%
{lMatem}) that%
\begin{equation}
A_{n}=H_{n}^{(1)}H_{n}^{(2)}\cdots H_{n}^{(m)}.  \label{Matexplic}
\end{equation}

Further, let $\mathbf{C}_{n}:=\left( C_{1}(n),\ldots ,C_{m}(n)\right)
^{\prime }$ be a random vector with components
\begin{equation*}
C_{i}(n):=\frac{dF_{n}^{(i)}(\mathbf{1;}\lambda )}{d\lambda }\Big|
_{\lambda =0} =\mathbf{E}_{\mathbf{F}_{n}}\varphi
_{i}(n),\,i=1,\ldots ,m
\end{equation*}%
and let $\mathbf{c}_{n}:=\left( c_{1}(n),\ldots ,c_{m}(n)\right) ^{\prime }$
be a random vector with components%
\begin{equation*}
c_{i}(n):=\frac{d\phi _{n}^{(i)}(\mathbf{1;}\lambda )}{d\lambda
}\Big| _{\lambda =0} =\mathbf{E}_{\mathbf{\phi }_{n}}\phi
_{i}(n),\,i=1,\ldots ,m.
\end{equation*}

Then, by (\ref{GG3}) $C_{m}(n)=c_{m}(n)$ and, for $i<m$%
\begin{equation*}
C_{i}(n)=c_{i}(n)+\sum_{k=i+1}^{m}h_{ik}(n)C_{k}(n).
\end{equation*}%
Hence we get $\mathbf{C}_{n}=\left( E-H_{n}^{\Delta }\right) ^{-1}\mathbf{c}%
_{n}$ where
\begin{equation*}
H_{n}^{\Delta }:=\left( h_{ij}(n)I\left( i<j\right) \right) _{i,j=1}^{m}
\end{equation*}%
is the upper triangular matrix generated by $H_{n}$.

Let, further, $\mathbf{L}_{n}:=\left( l_{ij}(n)\right) _{i,j=1}^{m}$ be a
random matrix with elements%
\begin{equation*}
l_{ij}(n)=\frac{\partial l_{n}^{(i)}(\mathbf{s})}{\partial
s_{j}}\Big| _{\mathbf{s}=\mathbf{1}} =\mathbf{E}_{l_{n}}\zeta
_{ij}(n),
\end{equation*}%
and let%
\begin{equation*}
p_{i}(n):=\frac{d\rho _{n}^{(i)}(\mathbf{1;}\lambda )}{d\lambda
}\Big| _{\lambda =0} =\mathbf{E}_{\mathbf{\rho }_{n}}\sigma
_{i}(n),\,i=1,\ldots ,m.
\end{equation*}%
Then, by (\ref{II2}) and (\ref{II3})%
\begin{equation*}
B_{j}(n)=\frac{\partial g_{n}(\mathbf{s})}{\partial s_{j}}\Big| _{%
\mathbf{s}=\mathbf{1}}
=\sum_{i=1}^{j-1}l_{ij}(n)+\sum_{i=1}^{m}\sum_{k=i+1}^{m}l_{ik}(n)a_{kj}(n)
\end{equation*}%
and
\begin{equation*}
D_{n}=\mathbf{E}_{G_{n}}\psi (n)=\frac{dG_{n}(\mathbf{1},\lambda
)}{d\lambda }\Big| _ {\lambda =0} =\left(
\sum_{i=1}^{m}p_{i}(n)+\sum_{k=i+1}^{m}l_{ik}(n)C_{k}(n)\right) .
\end{equation*}

The next two statements are easy consequences of Theorem \ref{TTmain} .

\begin{theorem}
\label{TTbusy}Assume that the MBPIFPRE associated with a BTPSFPRE
is subcritical and satisfies conditions of Theorem \ref{TTmain}.
Then there exist constants $C_{1},C_{2}\in (0,\infty )$ such that
\begin{equation}
C_{1}y^{-\kappa }\leq \mathbf{P}\left( \Theta _{P}>y\right) \leq
C_{2}y^{-\kappa },\,y\geq y_{0}.  \label{Busy2}
\end{equation}%
In particular, if the final product of any customer is its service time then
the tail distribution of the length $\Theta _{P}$ of a busy period of the
system satisfies (\ref{Busy2}).

\begin{corollary}
\label{TTpolSub}Under the conditions of Theorem \ref{TTbusy} $\mathbf{E}%
\Theta _{P}^{x}<\infty $ if and only if $x\in (0,\kappa )$.
\end{corollary}
\end{theorem}

\textbf{Proof of Theorem \ref{TTbusy}}. We associate with the
initial polling system an auxiliary polling system obtained from
the initial one by excluding certain customers and the final
product contributed by them. We split the exclusion process into
several stages.

At the first stage we exclude from the initial system all the
customers (and the final product contributed by them) which arrive
to the system during the switch-over times of the server. We call
these customers the customers of the first level.

At the second stage we exclude from the initial system all the
customers (and the final product contributed by them) which arrive
to the system during the service times of the first level
customers. We call these customers the customers of the second
level.

At the $k$-th stage we exclude from the initial system all the
customers (and the final product contributed by them) which arrive
to the system during the service times of the $(k-1)$-th level
customers.  And so on.

 Denote by  $\Phi $ -- the total amount of the final product accumulated in
the auxiliary branching type polling system. It is clear that
$\Theta _{P}\geq \Phi. $  As shown in \cite{Vat2009}, under the
conditions of Theorem \ref{TTbusy} $\mathbf{P}\left( \Phi
>y\right) \sim Ky^{-\kappa }$ as $y\rightarrow \infty $. Hence
\begin{equation*}
\mathbf{P}\left( \Theta _{P}>y\right) \geq \mathbf{P}\left( \Phi >y\right)
\geq C_{1}y^{-\kappa },y\geq y_{0}.
\end{equation*}%
This proves the needed bound from below.

To get the desired estimate from above note that $\Theta _{P}\leq
\Theta _{ext}$ where $\Theta _{ext}$ is the total amount of the
final product accumulated in the analyzed branching type polling
system within the generalized busy
period. By Theorem \ref{RResRep} $\Theta _{ext}\overset{d}{=}\Theta $ where $%
\Theta $ is the total amount of the final product produced in the
associated MBPIFPRE within a life period which starts at moment
$0$ by the birth of a type $1$ particle. Using now Theorem
\ref{TTmain} we obtain
\begin{equation*}
\mathbf{P}\left( \Theta _{P}>y\right) \leq \mathbf{P}\left( \Theta
>y\right) \sim C_{I}y^{-\kappa },\quad y\rightarrow \infty ,
\end{equation*}%
which gives the desired estimate from above.

For completeness we formulate the following statement concerning polling
systems with $\alpha >0$.

\begin{theorem}
\label{TTpolSuper}Assume that the MBPIFPRE associated with a BTPSFPRE is
such that its underlying MBPRE satisfies conditions of Theorem \ref{Ttan}
with $\alpha >0$ and, in addition,%
\begin{equation*}
\mathbf{P}\left( \min_{1\leq l\leq m}\mathbf{E}_{\mathbf{f}}\varphi
_{l}>0\right) >0.
\end{equation*}
If $\Theta _{P}$ is the total size of the final product accumulated in the
BTPSFPRE during a busy period then $\mathbf{P}\left( \Theta _{P}=\infty
\right) >0$. In particular, if the service time of any customer at any
station is positive with probability 1 then the busy period of the BTPSFPRE
is infinite with positive probability.
\end{theorem}

\textbf{Proof}. \ As shown in \cite{Vat2009} under the conditions of Theorem %
\ref{TTpolSuper} $\mathbf{P}\left( \Phi =\infty \right) >0$ where $\Phi $ is
the same as in the proof of Theorem \ref{TTbusy}. Hence the desired
statement follows.

\end{document}